\documentstyle[12pt,twoside]{article}
\font\teneufm=eufm10 scaled \magstep1
\font\seveneufm=eufm7 scaled \magstep1
\font\fiveeufm=eufm5  scaled \magstep1
\newfam\eufmfam
\textfont\eufmfam=\teneufm
\scriptfont\eufmfam=\seveneufm
\scriptscriptfont\eufmfam=\fiveeufm
\def\frak#1{{\fam\eufmfam\relax#1}}
\newfam\msbfam
\font\tenmsb=msbm10 scaled \magstep1  \textfont\msbfam=\tenmsb
\font\sevenmsb=msbm7 scaled \magstep1 \scriptfont\msbfam=\sevenmsb
\font\fivemsb=msbm5 scaled \magstep1  \scriptscriptfont\msbfam=\fivemsb
\def\Bbb{\fam\msbfam \tenmsb}

\def\RR{{\Bbb R}}
\def\CC{{\Bbb C}}
\def\QQ{{\Bbb Q}}
\def\NN{{\Bbb N}}
\def\ZZ{{\Bbb Z}}
\def\TT{{\Bbb T}}

\def\ra{\rightarrow}
\def\HollowBoxx #1#2#3{{\dimen0=#1 \advance\dimen0 by -#2
       \dimen1=#1 \advance\dimen1 by #3
        \vrule height 0pt depth #3 width #2
       \hskip -#3
       \vrule height #1 depth #3 width #3}}
 \def\LeftContraction{\mathord{\kern1.45pt \HollowBoxx{6pt}{3.5pt}{.4pt}}\,}
 \def\HollowBox #1#2#3{{\dimen0=#1 \advance\dimen0 by -#3
       \dimen1=#1 \advance\dimen1 by #3
        \vrule height #1 depth #3 width #3
        \vrule height 0pt depth #3 width #2
        \hskip -#3}}
 \def\RightContraction{\mathord{\, \HollowBox{6pt}{3.1pt}{.4pt}} \kern1.6pt}
\def\qed{{\hfill $\Box$}}
\newtheorem{theorem}{THEOREM}[section]
\newtheorem{corollary}[theorem]{Corollary}
\newtheorem{lemma}[theorem]{Lemma}

\newtheorem{proposition}[theorem]{Proposition}

\begin{document}

\begin{center}
{\Large \bf Proper Holomorphic Maps between
\medskip\\
Reinhardt Domains in ${\bf C}^2$}
\medskip \\
\normalsize A. V. Isaev and N. G. Kruzhilin

\end{center}

\begin{quotation} \small \sl We characterize pairs of bounded Reinhardt domains in $\CC^2$ between which there exists a proper holomorphic map and find all proper maps that are not elementary algebraic.
\end{quotation}

\pagestyle{myheadings}
\markboth{A. V. Isaev and N. G. Kruzhilin}{Proper Holomorphic Maps between Reinhardt Domains}

\setcounter{section}{-1}

\section{Introduction}
\setcounter{equation}{0}

Let $D_1$, $D_2$ be bounded Reinhardt domains in $\CC^2$ and $f:D_1\ra D_2$ a proper holomorphic map. Such maps are often {\it elementary algebraic}, that is, have the \lq\lq monomial\rq\rq\, form
$$
\begin{array}{lll}
z&\mapsto&\hbox{const}\,z^aw^b,\\ 
w&\mapsto&\hbox{const}\,z^cw^d,
\end{array}
$$
where $z$, $w$ denote variables in $\CC^2$, and $a,b,c,d$ are integers such that $ad-bc\ne 0$. For brevity we shall call such maps {\it elementary}\, maps. All elementary  maps are well-defined outside $I$, the union of the coordinate complex lines, but not necessarily at  points in $I$. The question of the existence of an elementary proper holomorphic map between two given domains is resolved by passing to the logarithmic diagrams of the domains. Several classes of domains between which only elementary proper holomorphic maps are possible have been described in \cite{Sp}. 

The aim of the present paper is to identify situations in which $f$ is not elementary and to explicitly describe all forms that the map $f$ and the domains $D_1$, $D_2$ may have in such cases. If $f$ is biholomorphic, then it can be represented as the composition of an elementary biholomorphism between $D_1$, $D_2$ and automorphisms of these domains (see \cite{Kr}, \cite{Sh}). Therefore, non-elementary biholomorphisms can occur only between domains equivalent by means of an elementary map and having non-elementary automorphisms, and hence are straightforward to determine. 

Proper maps that are not biholomorphic are harder to deal with. 
Non-elementary maps may occur, for example, if both $D_1$ and $D_2$ are bidiscs in which case at least one component of $f$ contains a Blaschke product with a zero away from the origin. In \cite{BP}, \cite{LS} the problem of describing non-elementary proper holomorphic maps was studied for complete Reinhardt domains, and it turns out that, apart from the example of bidiscs, such maps can only arise if $D_1$ and $D_2$ are certain pseudoellipsoids. On the other hand, all proper holomorphic maps between pseudoellipsoids in $\CC^n$ for $n\ge 2$ can be found using arguments from \cite{D-SP}. All proper holomorphic maps for another special class of domains (a generalization to higher dimensions of domains of the form (\ref{final1}) below) were determined in \cite{Lan}. Similarly to bidiscs, the only non-elementary maps for this class are expressed in terms of Blaschke products with at least one zero away from the origin.  

In this paper we describe all non-elementary proper holomorphic maps
between Reinhardt domains in $\CC^2$, as well as the corresponding pairs of domains. First of all, the map $f$ can be extended to a proper holomorphic map between the envelopes of holomorphy $\hat D_1$ and $\hat D_2$ of $D_1$ and $D_2$, respectively. Further, $f$ extends holomorphically to a neighborhood of $\partial \hat D_1\setminus I$. Next, one can  show that $f$ can be non-elementary only if $\partial \hat D_1\setminus I$ either consists of two or three
Levi-flat pieces or is a connected spherical hypersurface
(see Section \ref{preliminaries}).

The Levi-flat and spherical cases are considered in Sections \ref{leviflat} and \ref{sphericalcase}, respectively, and the results are summarized in Theorem \ref{mainresult} below. In the spherical case (see (iv)--(vi) of Theorem \ref{mainresult}) the map $f$ can be represented as the composition of three maps of special forms: two elementary maps and an automorphism of an intermediate Reinhardt domain. We note that a factorization result of a different kind for proper maps into the unit ball was obtained in \cite{KLS}. In the case when $D_1$ is a strongly pseudoconvex smoothly bounded Reinhardt domain in $\CC^n$ for $n\ge 2$ not intersecting the coordinate hyperplanes (while $D_2$ is not necessarily Reinhardt), another factorization theorem was proved in \cite{BD}. We also observe that the non-elementary proper holomorphic map between pseudoellipsoids in the example given in \cite{D-SP} factors as in (vi) of Theorem \ref{mainresult}.

Our results immediately imply that if there exists a proper holomorphic map between two bounded Reinhardt domains, then there also exists an elementary proper map between the domains (Corollary \ref{cor1}). Another consequence of our classification is that the domains described in (i)-(iii) of Theorem (\ref{mainresult}) are the only domains for which there exist non-elementary non-biholomorphic proper self-maps (Corollary \ref{cor2}). 

\begin{theorem}\label{mainresult}\sl Let $D_1$, $D_2$ be bounded Reinhardt domains in $\CC^2$ and $f: D_1\ra D_2$ a proper holomorphic map. Assume that $f$ is not elementary. Then one of the following holds:
\smallskip\\

\noindent (i) Up to permutation of the components of $f$ and the variables, the map $f$ has the form
\begin{equation}
\begin{array}{lll}
z&\mapsto&\hbox{const}\,z^a w^b B(A_1 z^{p_1}w^{q_1}),\\
w&\mapsto&\hbox{const}\,w^c,
\end{array}\label{finalform1}
\end{equation} 
where $a,b,c,p_1,q_1\in\ZZ$, $a>0$, $c>0$, $p_1>0$, $q_1\le 0$, $p_1$ and $q_1$ are relatively prime, $aq_1-bp_1\le 0$, and $B$ is a non-constant Blaschke product in the unit disc non-vanishing at 0. In this case $D_1$ either has the form 
\begin{equation}
\left\{(z,w)\in\CC^2: A_1|z|^{p_1}|w|^{q_1}<1,\, 0<|w|<C_1\right\},\label{final1}
\end{equation}
for some $C_1>0$, or is a bidisc (in the second case $b=0$, $p_1=1$, $q_1=0$ in (\ref{finalform1})). The domain $D_2$ is respectively either a domain of the form
$$
\left\{(z,w)\in\CC^2: A_2|z|^{p_2}|w|^{q_2}<1,\,0<|w|<C_2\right\},
$$
where $p_2,q_2\in\ZZ$ are relatively prime, $p_2>0$, $q_2\le 0$, $q_2/p_2=(aq_1-bp_1)/(cp_1)$, and $A_2>0$, $C_2>0$, or a bidisc.
\smallskip\\

\noindent (ii) Up to permutation of the components of $f$ and the variables, the map $f$ has the form (\ref{finalform1}), where $a,b,c,p_1,q_1\in\ZZ$, $a>0$, $c\ne 0$, $p_1>0$, $p_1$ and $q_1$ are relatively prime, $A_1>0$, and $B$ is a non-constant Blaschke product in the unit disc non-vanishing at 0. In this case the domains have the forms 
\begin{equation}
D_1=\left\{(z,w)\in\CC^2: A_1|z|^{p_1}|w|^{q_1}<1,\, E_1<|w|<C_1\right\},\label{final111}
\end{equation}
$$
D_2=\left\{(z,w)\in\CC^2: A_2|z|^{p_2}|w|^{q_2}<1,\,E_2<|w|<C_2\right\},
$$
where $p_2,q_2\in\ZZ$ are relatively prime, $p_2>0$, $q_2/p_2=(aq_1-bp_1)/(cp_1)$, and $C_1>0$, $E_1>0$, $A_2>0$, $C_2>0$, $E_2>0$.
\smallskip\\

\noindent (iii) Up to permutation of the components of $f$, the map $f$ has the form
$$
\begin{array}{lll}
z&\mapsto&\hbox{const}\,z^a B_1(Az),\\
w&\mapsto&\hbox{const}\, w^b B_2(Cw),
\end{array}
$$
where $a,b\in\ZZ$, $a\ge 0$, $b\ge 0$,  $A>0$, $C>0$, and $B_1$, $B_2$ are non-constant Blaschke products in the unit disc non-vanishing at 0. In this case $D_1$, $D_2$ are bidiscs. 
\smallskip\\

\noindent (iv) The map $f$ is a composition $f={\bf g}\circ{\bf f}\circ{\bf h}$, where ${\bf h}$ is an elementary map from $D_1$ into the domain
$
D:=\left\{(z,w)\in\CC^2: |w|>\exp\left(|z|^2\right)\right\},
$
${\bf f}$ is an automorphism of $D$, and ${\bf g}$ is an elementary map from a subdomain of $D$ onto $D_2$. Up to permutation of the variables, the map ${\bf h}$ has the form
$$
\begin{array}{lll}
z&\mapsto& \hbox{const}\, z^{a_1}w^{-b_1},\\
w&\mapsto& \hbox{const}\, w^{-c_1},
\end{array}
$$
where $a_1,b_1,c_1\in\NN$; the map ${\bf f}$ has the form
$$
\begin{array}{lll}
z&\mapsto& e^{it_1}z+s,\\
w&\mapsto&
e^{it_2}\exp\left(2\overline{s}e^{it_1}z+
|s|^2\right)w,
\end{array}
$$
where $t_1,t_2\in\RR$, $s\in\CC^*$; up to permutation of its components, the map ${\bf g}$ has the form
$$
\begin{array}{lll}
z&\mapsto& \hbox{const}\, z^{a_2}w^{-b_2},\\
w&\mapsto& \hbox{const}\, w^{-c_2},
\end{array}
$$
where $a_2, b_2,c_2\in\NN$. In this case the domains have the forms 
$$
\begin{array}{ll}
D_1=\displaystyle \Biggl\{(z,w)\in\CC^2: C_1'\exp\left(-E_1 |z|^{2a_1}|w|^{-2b_1}\right)&<|w|<\\
& \displaystyle C_1\exp\left(-E_1 |z|^{2a_1}|w|^{-2b_1}\right)\Biggr\},\\
\vspace{0mm}&\\
D_2=\displaystyle\Biggl\{(z,w)\in\CC^2: C_2'\exp\left(-E_2 |z|^{\frac{2}{a_2}}|w|^{-\frac{2b_2}{a_2c_2}}\right)&<|w|<\\
&\displaystyle C_2\exp\left(-E_2 |z|^{\frac{2}{a_2}}|w|^{-\frac{2b_2}{a_2c_2}}\right)\Biggr\},    
\end{array}    
$$
where $0\le C_1'<C_1$, $0\le C_2'<C_2$,  $E_1>0$, $E_2>0$.
\smallskip\\

\noindent (v) The map $f$ is a composition $f={\bf g}\circ{\bf f}\circ{\bf h}$, where ${\bf h}$ is an elementary map from $D_1$ into the domain
$\Omega^{\alpha}:=\bigl\{(z,w)\in\CC^2:|z|^2+|w|^{\alpha}<1\bigr\}$,
for some $\alpha>0$, ${\bf g}$ is an elementary map from a subdomain of $\Omega^{\alpha}$ onto $D_2$, and ${\bf f}$ is an automorphism of $\Omega^{\alpha}$. Up to permutation of the variables, the map ${\bf h}$ has the form
\begin{equation}
\begin{array}{lll}
z&\mapsto& \hbox{const}\, z^{a_1}w^{-b_1},\\
w&\mapsto& \hbox{const}\, w^{c_1},
\end{array}\label{mapppping1}
\end{equation}
where $a_1, b_1, c_1\in\ZZ$, $a_1>0$, $b_1\ge 0$, $c_1>0$; the map ${\bf f}$ has the form
$$
\begin{array}{lll}
z&\mapsto&\displaystyle e^{it_1}\frac{z-a}{1-\overline{a}z},\\
\vspace{0mm}&&\\
w&\mapsto&\displaystyle e^{it_2}\frac{(1-|a|^2)^{\frac{1}{\alpha}}}{(1-\overline{a}z)^{\frac{2}{\alpha}}}w,
\end{array}
$$
where $|a|<1$, $a\ne 0$, $t_1, t_2\in\RR$; up to permutation of its components, the map ${\bf g}$ has the form
\begin{equation}
\begin{array}{lll}
z&\mapsto& \hbox{const}\, z^{a_2}w^{b_2},\\
w&\mapsto& \hbox{const}\, w^{c_2},
\end{array}\label{mapppping2}
\end{equation}
where $a_2, b_2,c_2\in\ZZ$, $a_2>0$, $b_2\ge 0$, $c_2>0$.  In this case the domains have either the forms
$$
\begin{array}{l}
D_1=\Bigl\{(z,w)\in\CC^2: C_1|z|^{2a_1}+E_1 |w|^{\alpha c_1}<1,\Bigr\},\\
\vspace{0mm}\\
D_2=\left\{(z,w)\in\CC^2: C_2|z|^{\frac{2}{a_2}}+E_2 |w|^{\frac{\alpha}{ c_2}}<1\right\},
\end{array} 
$$
or the forms
$$
\begin{array}{ll}
D_1=\Biggl\{(z,w)\in\CC^2:& C_1|z|^{2a_1}|w|^{-2b_1}<1,\,E_1'\left(1-C_1|z|^{2a_1}|w|^{-2b_1}\right)^{\frac{1}{\alpha c_1}}<|w|<\\
&\displaystyle E_1\left(1-C_1|z|^{2a_1}|w|^{-2b_1}\right)^{\frac{1}{\alpha c_1}}\Biggr\},\\
\vspace{0mm}&\\
D_2=\Biggl\{(z,w)\in\CC^2:& C_2|z|^{\frac{2}{a_2}}|w|^{-\frac{2b_2}{a_2c_2}}<1,\,E_2'\left(1-C_2|z|^{\frac{2}{a_2}}|w|^{-\frac{2b_2}{a_2c_2}}\right)^{\frac{c_2}{\alpha}}<|w|<\\
&\displaystyle E_2\left(1-C_2|z|^{\frac{2}{a_2}}|w|^{-\frac{2b_2}{a_2c_2}}\right)^{\frac{c_2}{\alpha}}\Biggr\},
\end{array}
$$
for some $C_1>0$, $C_2>0$, $0\le E_1'<E_1$, $0\le E_2'<E_2$ (in the first case in (\ref{mapppping1}), (\ref{mapppping2}) we have $b_1=0$, $b_2=0$).
\smallskip\\

\noindent (vi) The map $f$ is a composition $f={\bf g}\circ{\bf f}\circ{\bf h}$, where ${\bf h}$ is an elementary map from $D_1$ onto the unit ball
$B^2:=\bigl\{(z,w)\in\CC^2:|z|^2+|w|^2<1\bigr\}$,
${\bf f}$ is an automorphism of $B^2$, and ${\bf g}$ is an elementary map from $B^2$ onto $D_2$. Up to permutation of the variables, the map ${\bf h}$ has the form
$$
\begin{array}{lll}
z&\mapsto& \hbox{const}\, z^{a_1},\\
w&\mapsto& \hbox{const}\, w^{b_1},
\end{array}
$$
where $a_1, b_1\in\NN$; the map ${\bf f}$ is such that ${\bf f}\left(B^2\cap{\cal L}_z\right)\not\subset B^2\cap I$, ${\bf f}\left(B^2\cap{\cal L}_w\right)\not\subset B^2\cap I$, where ${\cal L}_z:=\{z=0\}$, ${\cal L}_w:=\{w=0\}$, $I:={\cal L}_z\cup{\cal L}_w$; up to permutation of the variables, the map ${\bf g}$ has the form
$$
\begin{array}{lll}
z&\mapsto& \hbox{const}\, z^{a_2},\\
w&\mapsto& \hbox{const}\, w^{b_2},
\end{array}
$$
where $a_2, b_2\in\NN$.  In this case the domains have the forms
$$
\begin{array}{l}
D_1=\Bigl\{(z,w)\in\CC^2: C_1|z|^{2a_1}+E_1 |w|^{2b_1}<1\Bigr\},\\
\vspace{0mm}\\
D_2=\left\{(z,w)\in\CC^2: C_2|z|^{\frac{2}{a_2}}+E_2 |w|^{\frac{2}{b_2}}<1\right\},
\end{array}
$$
where $C_1>0$, $E_1>0$, $C_2>0$, $E_2>0$.
\end{theorem}

We will now state two corollaries of Theorem \ref{mainresult} mentioned earlier.

\begin{corollary}\label{cor1}\sl If $D_1$ and $D_1$ are bounded Reinhardt domains in $\CC^2$, and there exists a proper holomorphic map from $D_1$ onto $D_2$, then there also exists an elementary proper map from $D_1$ onto $D_2$.
\end{corollary}

\begin{corollary}\label{cor2}\sl Let $D$ be a bounded Reinhardt domain in $\CC^2$ that admits a non-elementary non-biholomorphic proper holomorphic self-map. Then $D$ either up to permutation of the variables has one of the forms (\ref{final1}), (\ref{final111}), or is a bidisc. 
\end{corollary}

On the other hand,  if a bounded pseudoconvex Reinhard
domain $D$ admits an elementary non-biholomorphic proper holomorphic self-map, then $D$ either is a bidisc, or up to permutation of the variables has one of the forms (\ref{final1}), (\ref{final111}), or the form
$$
\left\{(z,w)\in\CC^2: A|z|^{p}|w|^{q}<1, E<|z|^{p'}|w|^{q'}<C\right\},
$$
where $A>0$, $C>0$, $E\ge 0$, and $p,q,p'q'$ are integers satisfying conditions similar to those in (\ref{final1}), (\ref{final111}). This is easy to see from the observation that if the  logarithmic diagram of $D$ is unbounded, then the two asymptotes of its convex boundary define either the eigendirections of the linear part of the affine transformation of the logarithmic diagram corresponding to the elementary  map, or those of the square of this operator.

Before proceeding, we would like to acknowledge that this work started while the second author was visiting the Department of Mathematics, Australian National University.

\section{Preliminaries}\label{preliminaries}
\setcounter{equation}{0}

As we pointed out in the introduction, $f$ can be extended to a proper holomorphic map (that we also denote by $f$) between the envelopes of holomorphy $\hat D_1$, $\hat D_2$ of $D_1$, $D_2$, respectively (see\cite{Ke}). Further, it follows from \cite{B} (see also \cite{Lan}) that $f$ extends holomorphically to a neighborhood of $\partial \hat D_1\setminus I$, where $I:={\cal L}_z\cup{\cal L}_w$, ${\cal L}_z:=\{z=0\}$, ${\cal L}_w:=\{w=0\}$. Since $f$ is proper, it follows that $f(\partial \hat D_1\setminus I)\subset \partial \hat D_2$.

For every $p=(z,w)\in\CC^2$ let $\TT(p)$ be the torus $\{(e^{i\alpha}z,e^{i\beta}w)\in\CC^2:\alpha,\beta\in\RR\}$ and let $\TT$ be the standard torus $\TT((1,1))$. We shall think of $\TT$ as a group acting on $\CC^2$. Next, for $j=1,2$ we denote by $H_j$ the union of all locally holomorphically homogeneous real-analytic hypersurfaces lying in $\partial \hat D_j\setminus I$ and set $S_j:=\partial \hat D_j\setminus (H_j\cup I)$. Further, denote by $J_f$ the zero set of the Jacobian of $f$ in $\partial\hat D_1\setminus I$ and let $C_f:=f^{-1}\Bigl(f(\partial \hat D_1\setminus I)\cap I\Bigr)$. 
We shall start with the following lemma (see also \cite{LS} and \cite{Sp}).

\begin{lemma}\label{nonsmooth}\sl If $p\in S_1$, then $f(\TT(p))\subset \TT(f(p))$. In particular, $\TT(f(p))\not\subset I$.
\end{lemma}

\noindent {\bf Proof:} Assume that $f(\TT(p))\not\subset \TT(f(p))$ and let  $p'\in\TT(p)$ be a point close to $p$ such that $p'\not \in C_f\cup J_f$ and $f(\TT(p'))=f(\TT(p))$ is not tangent to $\TT(f(p'))$. Choose a neighborhood $U$ of $p'$ in which $f$ is biholomorphic and let $V:=f(U)$. We may assume that $f(\TT(q))\not\subset \TT(f(p'))$ for all $q\in U$. Let $T:=V\cap \TT(f(p'))$ and let $\gamma\subset f(\TT(p')) $ be the image of the orbit of $p'$  on $\TT(p')$ under the action of a 1-parameter subgroup of  $\TT$  such that $\gamma$ is not tangent to $T$. Consider now the set $\Gamma:=\cup_{s\in\gamma}\TT(s)$. This is clearly a real-analytic hypersurface in $\partial \hat D_2$, which is, moreover, locally holomorphically homogeneous
because we have on it actions of a 2-dimensional torus and a 1-parameter  group, and the orbits of one action are transversal to those of the other. Thus, $f(p')\in H_2$, and therefore  $p'\in  H_1$. This means that $p\in H_1$, which contradicts the assumptions of the lemma.
Hence $f(\TT(p))\subset \TT(f(p))$, as required.
\qed
\smallskip\\

By Lemma 4.4 of \cite{Sp} and Lemma \ref{nonsmooth}, if $S_1$ contains at least three distinct tori, the map $f$ is elementary. Therefore, from now on we assume that $S_1$ contains no more than two distinct tori.

Hypersurfaces making up $H_j$ are either strongly pseudoconvex or Levi flat, and we denote by $H_j^{\hbox{\tiny spher}}$, $H_j^{\hbox{\tiny non-spher}}$ and $H_j^{\hbox{\tiny flat}}$ the unions of all spherical (that is, locally biholomorphically equivalent to the unit sphere in $\CC^2$), non-spherical and Levi flat hypersurfaces from $H_j$, respectively, for $j=1,2$. 

We shall deal with the case $H_1^{\hbox{\tiny non-spher}}\ne\emptyset$ first. Note that locally holomorphically homogeneous non-spherical Reinhardt hypersurfaces do exist and Lemma 3.3 of \cite{Sp} stating otherwise is incorrect. Consider, for example, the non-spherical tube hypersurface
$$
{\bf T}:=\left\{(z,w)\in\CC^2:\hbox{Re}\, w=\left(\hbox{Re}\,z\right)^3,\,\,\hbox{Re}\,z>0\right\}.
$$
The base of ${\bf T}$ is an affinely homogeneous curve, and hence ${\bf T}$ is holomorphically homogeneous. The map $\Pi: (z,w)\mapsto(e^z,e^w)$ takes suitable portions of ${\bf T}$ to 
locally holomorphically homogeneous non-spherical Reinhardt hypersurfaces. 

As the following proposition shows, if $H_1^{\hbox{\tiny non-spher}}\ne\emptyset$, the map $f$ is elementary.       

\begin{proposition}\label{nonspher} \sl If $M$ is a locally holomorphically homogeneous strongly pseudoconvex non-spherical hypersurface in $H_1$, then $f(\TT(p))\subset\TT(f(p))$ for every $p\in M\setminus (C_f\cup J_f)$.
\end{proposition}

\noindent {\bf Proof:} A Reinhardt hypersurface $N\subset\CC^2\setminus I$ is locally biholomorphically equivalent to the tube hypersurface $T_N:=\log(N)+i\RR^2$ the base of which is the logarithmic diagram $\log(N)\subset\RR^2$ of $N$ (cf. the example above). If $N$ is real-analytic, strongly pseudoconvex, non-spherical and homogeneous, infinitesimal CR-transformations of $T_N$ form a 3-dimensional Lie algebra ${\frak g}_{T_N}$ (see, e.g., \cite{C}). Further, it follows from \cite{Lob} that the curve $\log(N)$ is locally affinely homogeneous. Taking into account that translations in the imaginary directions form a two-dimensional subalgebra ${\frak h}_{T_N}$ in ${\frak g}_{T_N}$, we see that ${\frak g}_{T_N}$ is generated by ${\frak h}_{T_N}$ and a one-dimensional algebra of local affine transformations of $\log(N)$. Hence ${\frak h}_{T_N}$ is an ideal in ${\frak g}_{T_N}$, and it is straightforward to observe that there are no other ideals in ${\frak g}_{T_N}$. Let ${\frak h}_N$ be the ideal corresponding to ${\frak h}_{T_N}$ in the Lie algebra ${\frak g}_N$ of all infinitesimal CR-transformations of $N$. The ideal ${\frak h}_N$ consists of infinitesimal transformations corresponding to the action of $\TT$ on $N$.

Fix now $p\in M\setminus (C_f\cup J_f)$. Clearly, $f$ maps a neighborhood of $p$ in $M$ biholomorphically onto a  hypersurface $M'\subset H_2^{\hbox{\tiny non-spher}}$. The homomorphism between ${\frak g}_{M}$ and ${\frak g}_{M'}$ induced by $f$ maps ${\frak h}_{M}$ into ${\frak h}_{M'}$. Hence we have $f(\TT(p))\subset\TT(f(p))$, as required.
\qed
\smallskip\\

It follows from Lemma 4.4 of \cite{Sp} and Proposition \ref{nonspher} that $f$ is elementary, if $H_1^{\hbox{\tiny non-spher}}\ne\emptyset$.  Thus from now on we shall assume that $H_1^{\hbox{\tiny non-spher}}=\emptyset$.

Assume now that $H_1^{\hbox{\tiny spher}}\ne\emptyset$ and let $M$ be a spherical hypersurface in $H_1$. It follows from the proof of Proposition 3.2 of \cite{Sp} that $f(\TT(p))\subset \TT(f(p))$ for all $p\in M$, if the closure of $M$ intersects $S_1$. In this case $f$ is again elementary by Lemma 4.4 of \cite{Sp}. 

To summarize, non-elementary proper holomorphic maps can only exist in the following two cases: either $H_1= H_1^{\hbox{\tiny flat}}$ and $S_1$ consists of one or two distinct tori, or $\partial \hat D_1\setminus I$ is a connected spherical hypersurface. These cases are considered in Sections \ref{leviflat} and \ref{sphericalcase}, respectively. 

\section{Levi Flat Case}\label{leviflat}
\setcounter{equation}{0}

In this section we assume that $H_1= H_1^{\hbox{\tiny flat}}$. Note that in this case the logarithmic diagram of $\hat D_1$ is an unbounded polygon with either one or two vertices, depending on the number (one or two) of tori in $S_1$.

{\bf The case of a single torus}. Let $f_1, f_2$ be the components of $f$ and assume first that  $S_1$ is a single torus $\TT_1$. Let $\log(\hat D_1)$ be the logarithmic diagram of $\hat D_1$. The set $H_1$ can be represented as the union of two distinct Levi flat hypersurfaces $L_1^1$, $L_1^2$ whose boundaries in $\CC^2\setminus I$ coincide with $\TT_1$. Further, since $\hat D_1$ is bounded, $\log(\hat D_1)$ is a sector lying in the interior of a right angle of the form $\{(x,y)\in\RR^2:x<x_0,\, y<y_0\}$ for some $(x_0, y_0)\in\RR^2$. We can describe it as follows
$$
\log(\hat D_1)=\left\{(x,y)\in\RR^2: a_1x-y>-\ln C, \, x+d_1y<-\ln A\right\},
$$
where $a_1\ge 0$, $d_1\le 0$, $a_1d_1> -1$ and $A,C>0$. Note that $\hat D_1$ can contain the origin only if $a_1=d_1=0$. 

Each of $L_1^1$ and $L_1^2$ is foliated by complex curves, and every such curve intersects $\TT_1$ along a real-analytic curve. Hence, we obtain two distinct families of curves ${\cal C}_1^j$, $j=1,2$, on $\TT_1$. If $\psi_1:\RR^2\ra\TT_1$ is the covering map, the inverse images of  ${\cal C}_1^j$ under $\psi_1$ are two distinct families of parallel lines ${\cal L}_1^j$ in $\RR^2$, $j=1,2$.

For $p\in\TT_1\setminus J_f$ consider the torus $\TT_2:=\TT(f(p))$. By Lemma \ref{nonsmooth} we obtain $f(\TT_1)\subset\TT_2$ and $\TT_2\not\subset I$. Clearly, if $U$ is a small neighborhood of $p$, then in a neighborhood of $f(p)$ the torus $\TT_2$ lies in the boundaries of two distinct Levi flat hypersurfaces $f(L_1^1\cap U)$ and $f(L_1^2\cap U)$. Hence $\TT_2$ entirely lies in the boundaries of two distinct Levi flat hypersurfaces $L_2^j$, $j=1,2$. The hypersurfaces $L_2^j$ produce two distinct families of curves ${\cal C}_2^j$ on $\TT_2$, and $f({\cal C}_1^j)\subset {\cal C}_2^j$, $j=1,2$. Each ${\cal C}_2^j$ is invariant under the action of $\TT$ on $\TT_2$, and therefore, if $\psi_2:\RR^2\ra\TT_2$ is the covering map, the inverse images of  ${\cal C}_2^j$ under $\psi_2$ are two distinct families of parallel lines ${\cal L}_2^j$ in $\RR^2$, $j=1,2$. Further, if $\tilde f=(\tilde f_1,\tilde f_2):\RR^2\ra\RR^2$ is the real-analytic covering map for $f|_{\TT_1}:\TT_1\ra\TT_2$, then $\tilde f({\cal L}_1^j)\subset {\cal L}_2^j$ for $j=1,2$.

Let $g$ be a linear transformation of $\RR^2$ mapping ${\cal L}_1^1$ and ${\cal L}_1^2$ into the families of horizontal and vertical lines, respectively, and let $h$ be a similar transformation for the families ${\cal L}_2^j$, $j=1,2$. Consider $\hat f=h\circ\tilde f\circ g^{-1}$. Clearly, $\hat f=(\hat f_1,\hat f_2)$ is a real-analytic map such that $\hat f_1$ is constant on every vertical line and $\hat f_2$ is constant on every horizontal line in $\RR^2$. Hence $\hat f_1$ is a function of $x$ and $\hat f_2$ is a function of $y$ alone. We choose $g$ to be the linear transformation with the matrix 
$$
\left(
\begin{array}{lr}
a_1 & -1\\
1 & d_1
\end{array}
\right),
$$
and let the matrix of $h$ be
$$
\left(
\begin{array}{ll}
a_2 & b_2\\
c_2 & d_2
\end{array}
\right).  
$$
Since $h\circ\tilde f=\hat f\circ g$, it follows that
$$
\begin{array}{lll}
a_2\tilde f_1(x)+b_2\tilde f_2(y)&=&\hat f_1(a_1x-y),\\
c_2\tilde f_1(x)+d_2\tilde f_2(y)&=&\hat f_2(x+d_1y).
\end{array}
$$
This implies that there exist holomorphic functions of one variable $F$ and $G$ such that in a neighborhood $U$ of $p\in\TT_1$ we have
\begin{equation}
\begin{array}{lll}
f_1^{a_2}f_2^{b_2}&=&F(z^{a_1}w^{-1}),\\
f_1^{c_2}f_2^{d_2}&=&G(zw^{d_1}).
\end{array}\label{main}
\end{equation}

We shall consider the case $a_1,d_1\in\QQ$ first. Let $a_1=a_1'/a''_1$, where $a_1'\ge 0$, $a''_1>0$ are relatively prime integers. For a fixed $\alpha_1\ne 0$ consider the curve $P_1^{\alpha_1}$ with the equation $z^{a_1'}w^{-a''_1}=\alpha_1$. The logarithmic diagram $\log(P_1^{\alpha_1})$ of $P_1^{\alpha_1}$ is a straight line parallel to one side of $\log(\hat D_1)$. We choose $\alpha_1$ such that $P_1^{\alpha_1}\cap \hat D_1\cap U\ne\emptyset$. The intersection $P_1^{\alpha_1}\cap \hat D_1$ is biholomorphically equivalent to either a disc or a punctured disc, and the equivalence is given by $\zeta_1\mapsto (\mu_1\zeta_1^{a''_1},\nu_1\zeta_1^{a_1'})$, where $\zeta_1$ is the variable in a disc of a suitable radius, and $\mu_1^{a_1'}\nu_1^{-a''_1}=\alpha_1$. We note that $P_1^{\alpha_1}\cap \hat D_1$ can be equivalent to a disc only if $a_1'=0$.      

We have $f_1^{a_2}f_2^{b_2}=\alpha_2:=F\left(\alpha_1^{1/a''_1}\right)$ on an open subset of $U\cap P_1^{\alpha_1}$ in which $f_1^{a_2}$ and $f_2^{b_2}$ are defined as single-valued holomorphic functions. Hence $f(P_1^{\alpha_1}\cap \hat D_1)$ is contained in $P_2^{\alpha_2}\cap \hat D_2$, where $P_2^{\alpha_2}$ is obtained by the analytic continuation of a connected component of the analytic set defined by the equation $z^{a_2}w^{b_2}=\alpha_2$ near $f(p)$. Since $P_1^{\alpha_1}\cap \hat D_1$ is closed, so is $f(P_1^{\alpha_1}\cap \hat D_1)$, and therefore $f(P_1^{\alpha_1}\cap \hat D_1)=P_2^{\alpha_2}\cap \hat D_2$, and $a_2$, $b_2$ are rationally dependent. Changing the function $F$ if necessary, we can assume that either $a_2\in\QQ$, $a_2\ge 0$ and $b_2=-1$, or $a_2=1$ and $b_2=0$.  Clearly, the restriction of $f$ to $P_1^{\alpha_1}\cap \hat D_1$ is proper. Furthermore, $P_2^{\alpha_2}\cap \hat D_2$ is equivalent to either a disc or a punctured disc. If $b_2=-1$ and $a_2=a_2'/a''_2$ for some relatively prime integers $a_2'\ge 0$, $a''_2>0$, then this equivalence has the form $\zeta_2\mapsto (\mu_2\zeta_2^{a''_2},\nu_2\zeta_2^{a_2'})$, where $\zeta_2$ is the variable in a disc of a suitable radius, and $\mu_2^{a_2'}\nu_2^{-a''_2}=\alpha_2^{a''_2}$. If $a_2=1$, $b_2=0$, the equivalence has the form $\zeta_2\mapsto (\alpha_2,\zeta_2)$.

Assume first that for some $\alpha_1$ such that $P_1^{\alpha_1}\cap \hat D_1\cap U\ne\emptyset$, the intersections $P_1^{\alpha_1}\cap \hat D_1$ and $P_2^{\alpha_2}\cap \hat D_2$ are equivalent to punctured discs $r_1\buildrel{\circ}\over{\Delta}$ and $r_2\buildrel{\circ}\over{\Delta}$, respectively (we denote by $\Delta$ and $\buildrel{\circ}\over{\Delta}$ the unit disc and the punctured unit disc, respectively). A proper holomorphic map between $r_1\buildrel{\circ}\over{\Delta}$ and $r_2\buildrel{\circ}\over{\Delta}$ has the form $\zeta_2= \hbox{const}\,\zeta_1^k$, where $k$ is a positive integer. Hence from the second equation in (\ref{main}) we obtain
$$
\hbox{const}\,\zeta_1^{\sigma}=G(\hbox{const}\,
\zeta_1^{\mu}),
$$
for all $\zeta_1$ in an open subset of $r_1\buildrel{\circ}\over{\Delta}$ and some non-zero $\sigma,\mu\in\RR$. This means that $G(t)=\hbox{const}\,t^{\tau}$ for some $\tau\in\RR$.

Assume next that for some $\alpha_1$ such that $P_1^{\alpha_1}\cap \hat D_1\cap U\ne\emptyset$, the intersections $P_1^{\alpha_1}\cap \hat D_1$ and $P_2^{\alpha_2}\cap \hat D_2$ are equivalent to discs $r_1\Delta$ and $r_2\Delta$, respectively (in this case $a_1'=0$, $a''_1=1$). A proper holomorphic map between $r_1\Delta$ and $r_2\Delta$ has the form $\zeta_2= r_2 B(\zeta_1/r_1)$, where $B$ is a Blaschke product in the unit disc. Hence from the second equation in (\ref{main}) we obtain
$$
\hbox{const}\,B(\zeta)^{\tau}=G(\hbox{const}\,
\zeta^{\mu}),
$$
for all $\zeta$ in an open subset of the unit disc $\Delta$ and some non-zero $\tau,\mu\in\RR$. This means that $G(t)=\hbox{const}\,B(\hbox{const}\,t^{\sigma})^{\tau}$ for some $\sigma,\tau\in\RR$.

In a similar way, considering the curves $Q_1^{\beta_1}$ and  $Q_2^{\beta_2}$ with the equations $z^{d''_1}w^{d_1'}=\beta_1$, where $d_1=d_1'/d''_1$ for some relatively prime integers $d_1'\le 0$, $d''_1>0$, and $z^{c_2}w^{d_2}=\beta_2:=G(\beta_1^{1/d''_1})$, respectively, we obtain that $F(t)=\hbox{const}\, t^{\rho}$, if $Q_1^{\beta_1}\cap \hat D_1$ and $Q_2^{\beta_2}\cap \hat D_2$ are equivalent to punctured discs, and $F(t)=\hbox{const}\,B(\hbox{const}\,t^{\eta})^{\rho}$ for some Blaschke product $B$ in the unit disc, if $Q_1^{\beta_1}\cap \hat D_1$ and $Q_2^{\beta_2}\cap \hat D_2$ are equivalent to discs (in the second case $d_1'=0$, $d''_1=1$), $\eta,\rho\in\RR$.

If $F(t)=\hbox{const}\, t^{\rho}$ and $G(t)=\hbox{const}\,t^{\tau}$, it follows from (\ref{main}) that $f$ is elementary. 

Let $F(t)=\hbox{const}\, t^{\rho}$ and $G(t)=\hbox{const}\,B(\hbox{const}\,t^{\sigma})^{\tau}$, where $B$ is a Blaschke product in the unit disc with a zero away from 0. In this case $a_1'=0$, $a''_1=1$. It now follows from (\ref{main}) that $f$ has either the form
\begin{equation}
\begin{array}{lll}
f_1(z,w)&=&\hbox{const}\,z^aw^b \tilde B(A^{d''_1} z^{d''_1}w^{d_1'}),\\
f_2(z,w)&=&\hbox{const}\,w^d,
\end{array}\label{form1}
\end{equation}
where $a,b,d\in\ZZ$ and $\tilde B$ is a non-constant Blaschke product in the unit disc non-vanishing at 0, or the form \begin{equation}
\begin{array}{lll}
f_1(z,w)&=&\hbox{const}\,z^aw^b \hat B(A^{d''_1} z^{d''_1}w^{d_1'}),\\
f_2(z,w)&=&\hbox{const}\,z^cw^d \tilde B(A^{d''_1} z^{d''_1}w^{d_1'}),
\end{array}\label{form2}
\end{equation}
where $a,b,c,d\in\ZZ$, $\tilde B$ is a non-constant Blaschke product in the unit disc non-vanishing at 0, $\hat B$ is either a Blaschke product in the unit disc with the same zeroes as $\tilde B$ or a constant, and $a$ may be nonzero only if $\hat B$ is non-constant. Forms (\ref{form1}) and (\ref{form2}) correspond to the cases $a_2=0$ and $a_2\ne 0$, respectively.  

We shall now show that form (\ref{form2}) can be simplified. Assume first that $\hat B$ is non-constant. Then $\hat D_2$ contains the origin, and $f^{-1}(0)$ contains the intersection of a curve of the form $z^{d''_1}w^{d_1'}=\hbox{const}$ with $\hat D_1$. Hence $f^{-1}(0)$ is not compact which contradicts the assumption that $f$ is proper. Thus, $\hat B\equiv \hbox{const}$, and therefore $a=0$. Hence (\ref{form2}) in fact only differs from (\ref{form1}) by permutation of the components of the maps. 

We shall now study form (\ref{form1}) of proper maps and the domains $\hat D_1$, $\hat D_2$ in more detail. For every $\alpha_1$ such that $|\alpha_1|>1/C$ we have $f(P_1^{\alpha_1}\cap \hat D_1)=P_2^{\alpha_2}\cap \hat D_2$, where $\alpha_2:=F\left(\alpha_1\right)$, and each of the curves $P_1^{\alpha_1}\cap \hat D_1$, $P_2^{\alpha_2}\cap \hat D_2$ is equivalent to either a disc or a punctured disc. However, $\tilde B$ has a zero away from 0, and therefore $P_2^{\alpha_2}\cap \hat D_2$ is in fact equivalent to a disc, and hence so is $P_1^{\alpha_1}\cap \hat D_1$. This shows that either 
\begin{equation}
\hat D_1=\left\{(z,w)\in\CC^2: A^{d''_1}|z|^{d''_1}|w|^{d_1'}<1,\, 0<|w|<C\right\},
\label{domain1}
\end{equation}
or
$$
\hat D_1=\left\{(z,w)\in\CC^2: |z|<1/A,\, |w|<C\right\}
$$
(in the second case we have $d_1'=0$, $d_1''=1$). Further, $f$ of the form (\ref{form1}) is a proper map from $\hat D_1$ onto a bounded Reinhardt domain only if $d>0$ and $ad_1'-bd_1''\le 0$.

It is straightforward to observe that there exists no proper subdomain of $\hat D_1$ mapped properly by $f$ onto a bounded Reinhardt domain and whose envelope of holomorphy coincides with $\hat D_1$. Thus, $D_1=\hat D_1$, and hence $D_2=\hat D_2$. If (\ref{domain1}) holds, then we have
$$
D_2=\left\{(z,w)\in\CC^2: \tilde A|z|^{\tilde d_1''}|w|^{\tilde d_1'}<1,\,0<|w|<\tilde C\right\},
$$
for some relatively prime integers $\tilde d_1',\tilde d_1''$, with $\tilde d_1'\le 0$, $\tilde d_1''>0$, such that $\tilde d_1'/\tilde d_1''=(ad_1'-bd_1'')/(dd_1'')$, and $\tilde A>0$, $\tilde C>0$. If $D_1$ is the bidisc, then $f$ can only be proper if $b=0$, and in this case $D_2$ is also a bidisc. We have thus obtained (i) of Theorem \ref{mainresult}.

The case $F(t)=\hbox{const}\,B(\hbox{const}\,t^{\eta})^{\rho}$, $G(t)=\hbox{const}\, t^{\tau}$, where $B$ is a Blaschke product in the unit disc with a zero away from 0, leads to the same description of $f$, $D_1$, $D_2$ up to permutation of the components of $f$ and the variables.

Let $F(t)=\hbox{const}\,B_1(\hbox{const}\,t^{\eta})^{\rho}$ and $G(t)=\hbox{const}\,B_2(\hbox{const}\,t^{\sigma})^{\tau}$, where $B_1$, $B_2$ are Blaschke products in the unit disc with zeroes away from 0. In this case $a_1'=d_1'=0$ and $a''_1=d''_1=1$. From (\ref{main}) we see that either $f$ has the form
\begin{equation}
\begin{array}{lll}
f_1(z,w)&=&\hbox{const}\,z^aw^b \tilde B_1(Az)\hat B_1(w/C),\\
f_2(z,w)&=&\hbox{const}\, w^d \hat B_2(w/C),
\end{array}\label{form3}
\end{equation}
where $a,b,d\in\ZZ$, $\tilde B_1$, $\hat B_2$ are non-constant Blaschke products in the unit disc non-vanishing at 0, $\hat B_1$ is either a Blaschke product in the unit disc with the same zeroes as $\hat B_2$ or a constant, and $b$ can be non-zero only if $\hat B_1$ is non-constant, or $f$ has the form
\begin{equation}
\begin{array}{lll}
f_1(z,w)&=&\hbox{const}\,z^aw^b \tilde B_1(Az)\hat B_1(w/C),\\
f_2(z,w)&=&\hbox{const}\, z^cw^d \tilde B_2 (Az)\hat B_2(w/C),
\end{array}\label{form4}
\end{equation}
where $a,b,c,d\in\ZZ$, $\hat B_1$ and $\tilde B_2$ are non-constant Blaschke products in the unit disc non-vanishing at 0, $\tilde B_1$ is either a Blaschke product in the unit disc with the same zeroes as $\tilde B_2$ or a constant, $\hat B_2$ is either a Blaschke product in the unit disc with the same zeroes as $\hat B_1$ or a constant, $a$ can be non-zero only if $\tilde B_1$ is non-constant, $d$ can be non-zero only if $\hat B_2$ is non-constant. Forms (\ref{form3}) and (\ref{form4}) correspond to the cases $a_2=0$ and $a_2\ne 0$, respectively.    

We shall now show that forms (\ref{form3}) and (\ref{form4}) can be simplified. Assume that in (\ref{form3}) $\hat B_1$ is non-constant. Then $\hat D_2$ contains the origin and $f^{-1}(0)$ contains the intersection of a complex line of the form $w=\hbox{const}$ with $\hat D_1$. Hence $f^{-1}(0)$ is not compact which contradicts the assumption that $f$ is proper. Therefore $\hat B_1\equiv\hbox{const}$, and hence $b=0$. A similar argument shows that in formula (\ref{form4}) we have $\tilde B_1\equiv\hbox{const}$ and $\hat B_2\equiv\hbox{const}$, and therefore $a=d=0$. Hence (\ref{form3}) reduces to
\begin{equation}
\begin{array}{lll}
f_1(z,w)&=&\hbox{const}\,z^a \tilde B_1(Az),\\
f_2(z,w)&=&\hbox{const}\, w^d \hat B_2(w/C),
\end{array}\label{form5}
\end{equation}
where $a,d\in\ZZ$, $a\ge 0$, $d\ge 0$, $\tilde B_1$, $\hat B_2$ are non-constant Blaschke products in the unit disc non-vanishing at 0, and (\ref{form4}) reduces to (\ref{form5}) up to permutation of the components of the maps.

Further, repeating the argument preceding formula (\ref{domain1}) we see that for every $\alpha_1$, $\beta_1$ with $|\alpha_1|>1/C$, $|\beta_1|<1/A$ the intersections $P_1^{\alpha_1}\cap \hat D_1$, $Q_1^{\beta_1}\cap \hat D_1$ are equivalent to discs, and therefore 
$$
\hat D_1=\left\{(z,w)\in\CC^2: |z|<1/A,\, |w|<C\right\}.
$$
Again, there exists no proper subdomain of $\hat D_1$ mapped properly by $f$ onto a bounded Reinhardt domain and whose envelope of holomorphy coincides with $\hat D_1$. Thus, $D_1=\hat D_1$, and hence $D_2=\hat D_2$. Therefore, $D_1$ and $D_2$ are bidiscs, and we have obtained (iii) of Theorem \ref{mainresult}.

Assume now that $a_1,d_1\not\in\QQ$. For a suitable $\alpha_1\ne 0$ consider the curve $P_1^{\alpha_1}$ obtained by the analytic continuation of the curve defined by the equation $z^{a_1}w^{-1}=\alpha_1$ in $U$. As before, we choose $\alpha_1$ to ensure that $P_1^{\alpha_1}\cap \hat D_1\cap U\ne\emptyset$. The intersection $P_1^{\alpha_1}\cap \hat D_1$ is not closed in $\hat D_1$ and is biholomorphically equivalent to to a half-plane; the equivalence is given by $\sigma_1:\zeta_1\mapsto (\exp(\zeta_1+\mu_1),\exp(a_1\zeta_1+\nu_1))$, where $\zeta_1$ is the variable in a suitable half-plane $R_1:=\left\{\zeta_1\in\CC:\hbox{Re}\,\zeta_1<s_1\right\}$, and $\exp(\mu_1a_1-\nu_1)=\alpha_1$.

As before, we observe that $f(P_1^{\alpha_1}\cap \hat D_1)$ lies in  $P_2^{\alpha_2}\cap \hat D_2$, where $P_2^{\alpha_2}$ for $\alpha_2:=F(\alpha_1)$ is obtained by the analytic continuation of a connected component of the set given by $z^{a_2}w^{b_2}=\alpha_2$ near $f(p)$. If $a_2$, $b_2$ were rationally dependent, the intersection $P_2^{\alpha_2}\cap \hat D_2$ would be closed in $\hat D_2$. Therefore $f^{-1}(P_2^{\alpha_2}\cap \hat D_2)$ would contain the closure of $P_1^{\alpha_1}\cap \hat D_1$ in $\hat D_1$ which is $|P_1^{\alpha_1}|\cap \hat D_1$, where $|P_1^{\alpha_1}|$ is the Levi flat hypersurface with the equation $|z|^{a_1}|w|^{-1}=|\alpha_1|$. Therefore, $a_2$, $b_2$ are in fact rationally independent, and $P_2^{\alpha_2}\cap \hat D_2$ is biholomorphically equivalent to either a half-plane or a strip, with the equivalence map $\sigma_2:\zeta_2\mapsto (\exp(-b_2\zeta_2+\mu_2),\exp(a_2\zeta_2+\nu_2))$, where $\zeta_2$ is the variable in either a suitable half-plane $R_2:=\left\{\zeta_2\in\CC:\hbox{Re}\,\zeta_2<s_2\right\}$, or a suitable strip $R_2':=\left\{\zeta_2\in\CC:s_2'<\hbox{Re}\,\zeta_2<s_2\right\}$, and $\exp(\mu_2a_2+\nu_2b_2)=\alpha_2$. Changing the function $F$ if necessary, we can assume that in the first case we have $a_2>0$, $b_2<0$. 

It is now straightforward to show that $f(P_1^{\alpha_1}\cap \hat D_1)=P_2^{\alpha_2}\cap \hat D_2$ and the restriction of $f$ to $P_1^{\alpha_1}\cap \hat D_1$ is proper. This restriction gives rise to a proper holomorphic map $\varphi:=\sigma_2^{-1}\circ f\circ \sigma_1$ either between $R_1$ and $R_2$, or between $R_1$ and $R_2'$. We shall now show that $\varphi$ is one-to-one. Assume the contrary and let $l_1$ be the line given by the equation $\hbox{Re}\,\zeta_1=s_1$. Since $\varphi$ is not one-to-one, $\varphi^{-1}(\infty)$ contains a point $\xi\in l_1$. Note that $\sigma_1^{-1}\Bigl(P_1^{\alpha_1}\cap(\partial \hat D_1\setminus I)\Bigr)=l_1$, and therefore $\sigma_1(\xi)\in\partial \hat D_1\setminus I$. In particular, $f$ is defined near $\sigma_1(\xi)$ and $f(\sigma_1(\xi))\in\partial \hat D_2$. On the other hand, consider in either $R_2$ or $R_2'$ a sequence $\{\xi_n\}$ converging to $\infty$, such that the sequence $\{\sigma_2(\xi_n)\}$ converges to a point in $\hat D_2$. Let $\{\xi_n'\}$ be a sequence in $R_1$ converging to $\xi$ such that $\varphi(\xi_n')=\xi_n$ for all $n$. Then $\{f(\sigma_1(\xi_n'))\}$ converges to a point in $\hat D_2$ which is impossible. Hence $\varphi$ is one-to-one. The above argument also shows that either $\varphi(l_1)=l_2$, or $\varphi(l_1)=l_2'$, where $l_2$, $l_2'$ are the lines given by the equations $\hbox{Re}\,\zeta_2=s_2$ and $\hbox{Re}\,\zeta_2=s_2'$, respectively. It follows that $P_2^{\alpha_2}\cap \hat D_2$ is in fact equivalent to $R_2$, $\varphi(l_1)=l_2$, and $\varphi(\zeta_1)=r\zeta_1 + q$, where $r>0$ and $q\in i\RR$.

Then from the second equation in (\ref{main}) we obtain
$$
\hbox{const}\,\exp(\sigma\zeta_1)=G\Bigl(\hbox{const}\,\exp(\mu\zeta_1)\Bigr)
$$
for all $\zeta_1$ in an open subset of $R_1$ and some non-zero $\sigma,\mu\in\RR$. Hence $G(t)=\hbox{const}\,t^{\tau}$ for some $\tau\in\RR$. Similarly, $F(t)=\hbox{const}\,t^{\rho}$ for some $\rho\in\RR$. It now follows from (\ref{main}) that $f$ is elementary.

We shall now assume that $a_1\not\in\QQ$ and $d_1\in\QQ$. Repeating the preceding arguments we obtain that $G(t)=\hbox{const}\,t^{\tau}$ for some $\tau\in\RR$ and either $F(t)=\hbox{const}\, t^{\rho}$ or $F(t)=\hbox{const}\,B(\hbox{const}\,t^{\eta})^{\rho}$ for some $\eta,\rho\in\RR$, where $B$ is a Blaschke product in the unit disc with a zero away from 0. In the first case we can show similarly to the above that $f$ is elementary. In the second case it is easy to see using (\ref{main}) that $f$ is necessarily a multi-valued map. This shows that the formula for $F$ does not actually contain a Blaschke product with a zero away from 0, and hence $f$ is elementary. Similarly, if $a_1\in\QQ$ and $d_1\not\in\QQ$, $f$ is elementary. 

{\bf The case of two tori}. Assume now that $S_1$ is a union of two tori. In this case $\log(\hat D_1)$ has the following form
$$
\begin{array}{lll}
\log(\hat D_1)&=&\Bigl\{(x,y)\in\RR^2: a_1x+b_1y>-\ln C, \, c_1x+d_1y<-\ln A, \\ &&u_1x+v_1y<-\ln E\Bigr\},
\end{array}
$$
for some $u_1$, $v_1$, where $a_1\ge 0$, $b_1\le 0$, $c_1\ge 0$, $d_1\le 0$, $b_1c_1\le a_1d_1$, and $A,C,E>0$. Note that $u_1$ and $v_1$ are not arbitrary: the line $u_1x+v_1y=0$ must intersect the other two \lq\lq to the left\rq\rq\, of their intersection point. 

The logarithmic diagram $\log(\hat D_1)$ has two vertices, and we shall concentrate on the one made by the lines $a_1x+b_1y=-\ln C$ and $u_1x+v_1y=-\ln E$ first. Let $\TT_1$ be the torus in $S_1$ corresponding to this vertex. As before, we can show that there exist holomorphic functions of one variable $F$ and $G$ such that in a neighborhood $U$ of $p\in\TT_1$
\begin{equation}
\begin{array}{lll}
f_1^{a_2}f_2^{b_2}&=&F(z^{a_1}w^{b_1}),\\
f_1^{u_2}f_2^{v_2}&=&G(z^{u_1}w^{v_1}),
\end{array}\label{main2}
\end{equation}
for some $a_2, b_2, u_2, v_2\in\RR$.

Assume first that both pairs $a_1$, $b_1$ and $u_1$, $v_1$ are rationally dependent. As before, we obtain that either $G(t)=\hbox{const}\, t^{\tau}$, or $G(t)=\hbox{const}\,B(\hbox{const}\,t^{\sigma})^{\tau}$, where $B$ is a Blaschke product in the unit disc with a zero away from 0, and $\sigma,\tau\in\RR$ (in the second case either $a_1=0$ or $b_1=0$). Similarly, considering the intersections  $Q_1^{\beta_1}\cap \hat D_1$, where $Q_1^{\beta_1}$ is the curve with the equation $z^{u_1}w^{v_1}=\beta_1$, we see $F(t)=\hbox{const}\,t^{\rho}$ for some $\rho\in\RR$. For the proof one must note that every connected component of $Q_1^{\beta_1}\cap \hat D_1$ is biholomorphically equivalent to an annulus with non-zero inner radius and that every proper map between two such annuli has the form $\zeta\mapsto\hbox{const}\,\zeta^k$, where $k\in\ZZ\setminus{0}$. 
 
For $G(t)=\hbox{const}\, t^{\tau}$ it follows from (\ref{main2}) that $f$ is elementary, and therefore we shall assume that $G(t)=\hbox{const}\,B(\hbox{const}\,t^{\sigma})^{\tau}$, where $B$ is a Blaschke product in the unit disc with a zero away from 0 (in this case either $a_1=0$ or $b_1=0$). Now (\ref{main2}) implies that up to permutation of its components, $f$ has either the form
\begin{equation}
\begin{array}{lll}
f_1(z,w)&=&\hbox{const}\,z^aw^b \tilde B(E^{u_1'/u_1} z^{u_1'}w^{v_1'}),\\
f_2(z,w)&=&\hbox{const}\,w^d,
\end{array}\label{form6}
\end{equation}    
where $a,b,d\in\ZZ$ and $\tilde B$ is a non-constant Blaschke product in the unit disc non-vanishing at 0, or the form
$$
\begin{array}{lll}
f_1(z,w)&=&\hbox{const}\,z^aw^b \tilde B(E^{v_1'/v_1}z^{u_1'}w^{v_1'}),\\
f_2(z,w)&=&\hbox{const}\,z^c,
\end{array}
$$  
where $a,b,c\in\ZZ$ and $\tilde B$ is a non-constant Blaschke product in the unit disc non-vanishing at 0. These forms correspond to the cases $a_1=0$ and $b_1=0$, respectively. In the first case $u_1> 0$, and $u_1'>0$, $v_1'$ are relatively prime integers such that $v_1/u_1=v_1'/u_1'$. In the second case $v_1>0$, and $u_1'$, $v_1'>0$ are relatively prime integers such that $u_1/v_1=u_1'/v_1'$. The above forms are obtained from one another by permutation of the variables, and we shall assume that (\ref{form6}) holds.

For $a_1=0$ the image of $\hat D_1$ under a map of the form (\ref{form6}) is a Reinhardt domain only if $c_1=0$, and we obtain
$$
\hat D_1=\left\{(z,w)\in\CC^2: E^{u_1'/u_1} |z|^{u_1'}|w|^{v_1'}<1,\, A^{-1/d_1}<|w|<C^{-1/b_1}\right\},
$$        
A map of the form (\ref{form6}) is a proper map from $\hat D_1$ onto a Reinhardt domain only if $d\ne 0$ and $a\ge 0$. As before, there exists no proper subdomain of $\hat D_1$ mapped properly by $f$ onto a bounded Reinhardt domain and whose envelope of holomorphy coincides with $\hat D_1$. Thus, $D_1=\hat D_1$, and hence $D_2=\hat D_2$. Then we have
$$
D_2=\left\{(z,w)\in\CC^2: \tilde E|z|^{\tilde u_1'}|w|^{\tilde v_1'}<1,\,\tilde A<|w|<\tilde C\right\},
$$
for some relatively prime integers $\tilde u_1',\tilde v_1'$, with $\tilde u_1'> 0$, such that $\tilde v_1'/\tilde u_1'=(av_1'-bu_1')/(du_1')$, and  $\tilde A>0$, $\tilde C>0$, $\tilde E>0$. We have thus obtained (ii) of Theorem \ref{mainresult}.

Assume now that $a_1$, $b_1$ are rationally dependent and $u_1$, $v_1$ are rationally independent. Then, as before, either $G(t)=\hbox{const}\, t^{\tau}$, or $G(t)=\hbox{const}\,B(\hbox{const}\,t^{\sigma})^{\tau}$, where $B$ is a Blaschke product in the unit disc with a zero away from 0, and $\sigma,\tau\in\RR$. Considering the intersections  $Q_1^{\beta_1}\cap \hat D_1$, where $Q_1^{\beta_1}$ is the curve with the equation $z^{u_1}w^{v_1}=\beta_1$, we see that $F(t)=\hbox{const}\,t^{\rho}$ for some $\rho\in\RR$. For the proof one must note that every connected component of $Q_1^{\beta_1}\cap \hat D_1$ is equivalent to a strip, with the equivalence map of the form $\zeta\mapsto (\exp(-v_1\zeta+\mu_1),\exp(u_1\zeta+\nu_1))$ with $\exp(\mu_1u_1+\nu_1v_1)=\beta_1$, and every proper map between two strips has the form $\zeta\mapsto r\zeta + q$, where $r\ne 0$ and $q\in i\RR$. If $G(t)=\hbox{const}\,t^{\tau}$ for some $\tau\in\RR$, then it follows from (\ref{main2}) that $f$ is elementary. If $G(t)=\hbox{const}\,B(\hbox{const}\,t^{\sigma})^{\tau}$, where $B$ is a Blaschke product in the unit disc with a zero away from 0, then it is easy to see from (\ref{main2}) that $f$ is necessarily a multi-valued map. This shows that the formula for $G$ does not actually contain a Blaschke product with a zero away from 0, and hence $f$ is elementary. 

If $a_1$, $b_1$ are rationally independent, we obtain $F(t)=\hbox{const}\,t^{\rho}$, $G(t)=\hbox{const}\,t^{\tau}$, with $\rho,\tau\in\RR$. In this case (\ref{main2}) shows that $f$ is elementary.

\section{Spherical Case}\label{sphericalcase}
\setcounter{equation}{0}

Assume now that $U_1:=\partial \hat D_1\setminus I$ is a connected real-analytic spherical hypersurface. Consider the logarithmic diagram $\log(\hat D_1)$ of $\hat D_1$ and the tube domain $T_1\subset\CC^2$ with base $\log(\hat D_1)\subset\RR^2$, that is, $T_1=\log(\hat D_1)+i\RR^2$. The domain $T_1$ covers $\hat D_1\setminus I$ by means of the map $\Pi:(z,w)\mapsto (e^z,e^w)$. Clearly, for every $p\in \hat D_1\setminus I$, the fiber $\Pi^{-1}(p)$ is preserved by the Abelian group $G$ of translations of $\CC^2$ of the form $(z,w)\mapsto (z+i2\pi  n, w+i2\pi  m)$, $n,m\in\ZZ$. The group $G$ has two generators: $(z,w)\mapsto (z+i2\pi, w)$ and $(z,w)\mapsto (z, w+i2\pi)$. We denote these maps by $\Lambda^z$ and $\Lambda^w$, respectively.  

Since $L_1:=\partial T_1$ is a closed spherical tube hypersurface, it follows from \cite{DY} that there exists an affine transformation $F_1$ of the form
\begin{equation}
\left(
\begin{array}{l}
z\\
w
\end{array}
\right)\mapsto
A
\left(
\begin{array}{l}
z\\
w
\end{array}
\right)+b,\label{affine}
\end{equation}
where $A\in GL_2(\RR)$, $b\in\RR^2$, that maps $L_1$ onto one of the four hypersurfaces defined by the following equations
\begin{equation}
\begin{array}{ll}
(1)& \hbox{Re}\, w=(\hbox{Re}\,z)^2,\\
(2)& \hbox{Re}\, w=\exp(2\hbox{Re}\,z),\\
(3)& \cos\left(\hbox{Re}\, w\right)=\exp(\hbox{Re}\,z),\\
(4)& \exp(2\hbox{Re}\,z)+\exp(2\hbox{Re}\,w)=1.
\end{array}\label{fourtypes}
\end{equation}
Let $\tilde T_1:=F_1(T_1)$. Clearly, $\tilde T_1$ covers $\hat D_1\setminus I$ by means of the map $\Pi_1:=\Pi\circ F_1^{-1}$, and for each $p\in \hat D_1\setminus I$ the group $G_1:=F_1\circ G \circ F_1^{-1}$ preserves the fiber $\Pi_1^{-1}(p)$.

Further, for each hypersurface listed in (\ref{fourtypes}) one can explicitly write a corresponding locally biholomorphic map onto a portion of the unit sphere $S^3$ respectively as follows (see \cite{DY})
\begin{equation}
\begin{array}{lll}
(1)& \displaystyle z\mapsto \frac{z}{\sqrt{2}},&\displaystyle w\mapsto w-\frac{z^2}{2},\\
\vspace{0mm}&\\
(2)& \displaystyle z\mapsto e^z,&\displaystyle w\mapsto w,\\
\vspace{0mm}&\\
(3)& \displaystyle z\mapsto \exp\left(\frac{z+iw}{2}\right),&
\displaystyle w\mapsto\exp(iw),\\
\vspace{0mm}&\\
(4)& \displaystyle z\mapsto \frac{e^z}{e^w-1},&
\displaystyle w\mapsto -\frac{e^w+1}{e^w-1}. 
\end{array}\label{fourmaps}
\end{equation}
In formulas (\ref{fourmaps}), $S^3$ punctured at a point is realized as the hypersurface with the equation $\hbox{Re}\, w=|z|^2$. The deleted point in this realization is at infinity and we denote it by $p_{\infty}$.

Let $\tilde L_1:=\partial\tilde T_1$ and let $\theta_1$ be the map from list (\ref{fourmaps}) corresponding to $\tilde L_1$. In case (1) $\tilde L_1$ is mapped by $\theta_1$ onto $S^3\setminus\{p_{\infty}\}$, in cases (2), (3) onto $S^3\setminus\left(\{p_{\infty}\}\cup{\cal L}_z\right)$, and in case (4) onto $S^3\setminus\left(\{p_{\infty}\}\cup{\cal L}_z\cup\{w=1\}\right)$. We also point out that in case (1) the map $\theta_1$ takes $\tilde T_1$ onto the unit ball $B^2$ realized as $\{(z,w)\in\CC^2:\hbox{Re}\, w> |z|^2\}$, in cases (2), (3) onto $B^2\setminus{\cal L}_z$, and in case (4) onto $B^2\setminus\left({\cal L}_z\cup\{w=1\}\right)$.  

For $h\in G_1$ consider now the locally defined map $\theta_1\circ h \circ \theta_1^{-1}$ from $S^3$ into itself. It extends to an automorphism $\hat h$ of $B^2$, and hence $G_1$ gives rise to a subgroup $\hat G_1$ of the group $\hbox{Aut}(B^2)$ of holomorphic automorphisms of $B^2$. The group $\hat G_1$ is clearly Abelian and has at most two generators.

Formulas (\ref{fourmaps}) yield the following descriptions of transformations in the group $\hat G_1$ in each of the four cases. In the expressions below the vector $(\alpha_1,\alpha_2)$ varies over a lattice in $\RR^2$.

\begin{equation}
\begin{array}{lll}
(1) &z\mapsto \displaystyle z+i\alpha_1,& \displaystyle w\mapsto -2i\alpha_1z+w+\alpha_1^2+i\alpha_2,\\
\vspace{0mm}&\\
(2)& z\mapsto \displaystyle e^{i\alpha_1}z,& \displaystyle w\mapsto w+i\alpha_2,\\
\vspace{0mm}&\\
(3)& z\mapsto \displaystyle e^{i\alpha_1+\alpha_2}z,& \displaystyle w\mapsto e^{2\alpha_2}w,\\
\vspace{0mm}&\\
(4)& z\mapsto\displaystyle\frac{2e^{i\alpha_1}z}{1+e^{i\alpha_2}+(1-e^{i\alpha_2})w},
& \displaystyle w\mapsto\displaystyle\frac{(e^{i\alpha_2}+1)w+1-e^{i\alpha_2}}{1+e^{i\alpha_2}+(1-e^{i\alpha_2})w}.
\end{array}\label{group}
\end{equation}
  
Next, since $f$ is locally biholomorphic at the points of $U_1\setminus J_f$, the set $H_2$ contains the real-analytic spherical hypersurface $U_2:=f(U_1\setminus(C_f\cup J_f))$. Let $T_2$ be the covering tube domain for $\hat D_2\setminus I$ and $L_2:=\Pi^{-1}(U_2)$ the portion of $\partial T_2$ covering $U_2$. By \cite{DY} there is an affine transformation $F_2$ of the form (\ref{affine}) mapping $L_2$ onto an open tube subset of one of hypersurfaces (\ref{fourtypes}). Let $\tilde L_2:=F_2(L_2)$. Clearly, $\tilde L_2$ covers $U_2$ by means of the map $\Pi_2:=\Pi\circ F_2^{-1}$, and for every $p\in U_2$ the group $G_2:=F_2\circ G \circ F_2^{-1}$ preserves the fiber $\Pi_2^{-1}(p)$. Let $\theta_2$ be the map from list (\ref{fourmaps}) corresponding to $\tilde L_2$. As for the group $G_1$, from every $h\in G_2$ by using the map $\theta_2$ we can produce $\hat h\in\hbox{Aut}(B^2)$, and therefore $G_2$ gives rise to a Abelian subgroup $\hat G_2\subset \hbox{Aut}(B^2)$ with at most two generators. In each of the four cases $\hat G_2$ is described by formulas (\ref{group}).

The map $f$ induces a homomorphism from $\hat G_1$ into $\hat G_2$ as follows. We fix $p_1\in U_1\setminus(C_f\cup J_f)$ and let  $p_2:=f(p_1)$. Clearly, $p_2\in U_2$. For $g_1\in G_1$ we choose $p_1', p_1''\in \Pi_1^{-1}(p_1)$ such that $g_1(p_1')=p_1''$ (note that $g_1$ is fully determined by this condition). Now fix a curve $\tilde\gamma_1\subset\Pi_1^{-1}(U_1\setminus (C_f\cup J_f))$ from $p_1'$ to $p_1''$ and let $\gamma_1:=\Pi_1(\tilde\gamma_1)$. Clearly, $\gamma_1$ is a closed curve in $U_1\setminus (C_f\cup J_f)$ passing through $p_1$ and $\gamma_2:=f(\gamma_1)$ is a closed curve in $U_2$ passing through $p_2$. For a fixed $p_2'\in\Pi_2^{-1}(p_2)$ we now consider a curve $\tilde\gamma_2\subset\tilde L_2$ originating at $p_2'$ such that $\Pi_2(\tilde\gamma_2)=\gamma_2$. Let $p_2''\in\Pi_2^{-1}(p_2)$ be the other endpoint of $\tilde\gamma_2$ and $g_2\in G_2$ be the map such that $g_2(p_2')=p_2''$. The correspondence $g_1\mapsto g_2$ defines a map from $G_1$ into $G_2$ which induces a map $\Phi: \hat G_1\ra \hat G_2$, $\Phi(\hat g_1)=\hat g_2$ (we show in the next paragraph that $\Phi$ is indeed a well-defined map). 

Denote by $\psi$ an analytic element of $\theta_1^{-1}$ defined near $\hat p_1:=\theta_1(p_1')$ such that $\psi(\hat p_1)=p_1'$, and by $\pi$ an analytic element of $\Pi_2^{-1}$ defined near $p_2$ such that $\pi(p_2)=p_2'$. Consider the map $\theta_2 \circ \pi \circ f \circ \Pi_1 \circ \psi$ mapping biholomorphically a neighborhood of $\hat p_1$ in $S^3$ onto a neighborhood of $\hat p_2:=\theta_2(p_2')$ in $S^3$. This map extends to an automorphism $\varphi$ of $B^2$, and one immediately observes that $\Phi(\hat g_1)= 
\varphi \circ \hat g_1 \circ \varphi^{-1}$ for all $\hat g_1\in\hat G_1$. This shows, in particular, that $\Phi$ is independent of the choice of the curves $\tilde \gamma_1$ and is single-valued. Clearly, $\Phi$ is a homomorphism. 

We now require the following result.

\begin{lemma}\label{finiteindex} \sl $\Phi(\hat G_1)$ is a finite-index subgroup of $\hat G_2$.
\end{lemma}

\noindent {\bf Proof:} Since $f$ is proper, $f^{-1}(p_2)$ consists of finitely many points, say, $p_1,q^1,\dots,q^k\in U_1$, $k\ge 0$. Let $\Gamma_1^1,\dots,\Gamma_1^k$ be curves in $U_1\setminus (C_f\cup J_f)$ joining respectively $q^1,\dots, q^k$ with $p_1$. Clearly, $\Gamma_2^j:=f(\Gamma_1^j)$, $j=1,\dots, k$, are closed curves in $U_2$ passing through $p_2$. As before, each curve $\Gamma_2^j$ gives rise to an element $g_2^j$ of $G_2$ and, consequently, to an element $\hat g_2^j$ of $\hat G_2$. 

Fix $g_2\in G_2$ and let $p_2'':=g_2(p_2')$. Let $\tilde\Gamma_2\subset \tilde L_2$ be a curve from $p_2'$ to $p_2''$ and let $\Gamma_2:=\Pi_2(\tilde\Gamma_2)$. Clearly, $
\Gamma_2$ is a closed curve in $U_2$ passing through $p_2$. Consider the curve $\Gamma_1\subset U_1\setminus (C_f\cup J_f)$ originating at $p_1$ such that $f(\Gamma_1)=\Gamma_2$.

If $\Gamma_1$ is closed, it gives rise to an element $\hat g_1$ of $\hat G_1$, and we obviously have $\hat g_2=\Phi(\hat g_1)$. Hence $\hat g_2\in\Phi(\hat G_1)$ in this case.

Assume now that $\Gamma_1$ is not closed and let $q^s$, for some $1\le s\le k$, be its other endpoint. Let $g_1$ be the element of $G_1$ corresponding to the closed curve obtained by joining $\Gamma_1$ and $\Gamma_1^s$. Then we clearly have $\Phi(\hat g_1)=\hat g_2+\hat g_2^s$, and hence $\hat g_2\in-\hat g_2^s+\Phi(\hat G_1)$.

We have thus shown that for $\hat g_2\in\hat G_2$ we either have $\hat g_2\in\Phi(\hat G_1)$, or $\hat g_2\in-\hat g_2^s+\Phi(\hat G_1)$ for some $1\le s\le k$. Therefore, $\Phi(\hat G_1)$ is of finite index in $\hat G_2$.

The proof of the lemma is complete.\qed
\smallskip\\

The following result imposes constraints on the possible forms of $\tilde L_1$ and $\tilde L_2$.

\begin{proposition}\label{nottwofour}\sl We have $\tilde L_2\subset \tilde L_1$, and the map $f$ is elementary unless $\tilde L_1$ is either hypersurface (2) or hypersurface (4) of (\ref{fourtypes}).
\end{proposition}

\noindent {\bf Proof:} First of all, we shall prove the first assertion. Assume that $\tilde L_2\not \subset \tilde L_1$. Then we shall show that the group $\hat G_1$ either cannot be conjugate in $\hbox{Aut}(B^2)$ to a subgroup of $\hat G_2$ or can only be conjugate to a subgroup of $\hat G_2$  of infinite index. It will therefore follow from Lemma \ref{finiteindex} that there exists no proper holomorphic map from $\hat D_1$ onto $\hat D_2$, contradicting our assumptions. Below we consider all possibilities for $\hat G_1$, $\hat G_2$ (see (\ref{group})).

Assume that $\tilde L_1$ is hypersurface (1) and $\tilde L_2$ lies in hypersurface (2). In this case the only fixed point of each of $\hat G_1$ and $\hat G_2$ in $\overline{B^2}$ is the point $p_{\infty}\in S^3$ at infinity. If $\varphi\circ\hat G_1\circ\varphi^{-1}\subset\hat G_2$ for some $\varphi\in\hbox{Aut}(B^2)$, then $\varphi(p_{\infty})=p_{\infty}$, i.e. $\varphi$ is affine. The general form of affine automorphisms of $B^2$ is as follows
\begin{equation}
\begin{array}{lll}
z&\mapsto&\lambda e^{it}z+\zeta,\\
\vspace{0mm}&&\\
w&\mapsto&\lambda^2w+2\overline{\zeta}\lambda e^{it}z+|\zeta|^2+i\mu,
\end{array}\label{formphi}
\end{equation}  
where $\lambda>0$, $\zeta\in\CC$, $t,\mu\in\RR$. It is now straightforward to show that $\hat G_1$ cannot be conjugate to a subgroup of $\hat G_2$ by means of an automorphism of the form (\ref{formphi}). The same argument works for the case when $\tilde L_1$ is hypersurface (2) and $\tilde L_2$ lies in hypersurface (1). 

Further, if $\tilde L_1$ is one of hypersurfaces (1), (2) and $\tilde L_2$ lies in hypersurface (3), the group $\hat G_1$ cannot be conjugate to a subgroup of $\hat G_2$ since $\hat G_1$ has only one fixed point in $S^3$ (the point $p_{\infty}$), whereas $\hat G_2$ has two (0 and $p_{\infty}$).

Let $\tilde L_1$  be hypersurface (3) and assume that $\tilde L_2$ lies in one of hypersurfaces (1), (2). Then $\hat G_1$ has two fixed points in $S^3$ (0 and $p_{\infty}$), and the only fixed point of $\hat G_2$ is $p_{\infty}$. It is clear from formula (\ref{group}) that $\hat G_2$ contains non-trivial elements fixing a point in $S^3$ other than $p_{\infty}$ only if $\tilde L_2$ lies in hypersurface (2); for such elements $\alpha_2=0$. However, a subgroup of $\hat G_2$ containing only elements satisfying this condition has infinite index in $\hat G_2$. Hence, $\hat G_1$ cannot be conjugate to a finite-index subgroup of $\hat G_2$.

Next, if $\tilde L_1$ is one of hypersurfaces (1), (2), (3) and $\tilde L_2$ lies in hypersurface (4), the group $\hat G_1$ cannot be conjugate to a subgroup of $\hat G_2$ since $\hat G_1$ does not have any fixed points in $B^2$, whereas $\hat G_2$ fixes the point $(0,1)\in B^2$.

Finally, let $\tilde L_1$ be hypersurface (4) and assume that $\tilde L_2$ lies in one of hypersurfaces (1), (2), (3). Then $\hat G_1$ has a fixed point in $B^2$ and $\hat G_2$ fixes no point in $B^2$. It is clear from formula (\ref{group}) that $\hat G_2$ contains non-trivial elements fixing a point in $B^2$ only if $\tilde L_2$ lies in either hypersurface (2) or hypersurface (3); for such elements $\alpha_2=0$. However, a subgroup of $\hat G_2$ containing only elements  satisfying this condition has infinite index in $\hat G_2$. Hence, $\hat G_1$ cannot be conjugate to a finite-index subgroup of $\hat G_2$.

We thus have shown that $\tilde L_2\subset \tilde L_1$. We shall now consider the two possibilities for $\tilde L_1$. In what follows we denote the map $\theta_1=\theta_2$ by $\theta$. 

Let $\tilde L_1$ be hypersurface (1) and $\Lambda^z_j:=\theta\circ F_j\circ\Lambda^z\circ F_j^{-1}\circ\theta^{-1}$, $\Lambda^w_j:=\theta\circ F_j\circ\Lambda^w\circ F_j^{-1}\circ\theta^{-1}$ be generators of $\hat G_j$, $j=1,2$. Since $\Phi(\hat G_1)\subset\hat G_2$, it follows that
$$
\begin{array}{l}
\varphi\circ\Lambda^z_1\circ\varphi^{-1}=(\Lambda^z_2)^{a_1}\circ(\Lambda^w_2)^{a_2},\\
\varphi\circ\Lambda^w_1\circ\varphi^{-1}=(\Lambda^z_2)^{b_1}\circ(\Lambda^w_2)^{b_2},
\end{array}
$$
for some $a_1,a_2,b_1,b_2\in\ZZ$, such that $a_1b_2-a_2b_1\ne 0$.

Consider the maps $\tilde\Lambda^z_1:=\theta^{-1}\circ\varphi\circ\Lambda^z_1\circ\varphi^{-1}\circ\theta$, $\tilde\Lambda^w_1:=\theta^{-1}\circ\varphi\circ\Lambda^w_1\circ\varphi^{-1}\circ\theta$. They generate a subgroup of the group $G_2$. Let $F$ be the linear map such that $F\circ\tilde\Lambda^z_1\circ F^{-1}=\Lambda^z$ and $F\circ\tilde\Lambda^w_1\circ F^{-1}=\Lambda^w$. 

We shall now introduce an intermediate domain $D$ through which the map $f$ can be factored. Let $T:=F(\tilde T_1)$ and $D:=\Pi(T)$. Clearly, $D$ is a Reinhardt domain. We define a biholomorphic map ${\bf f}$ from $\hat D_1\setminus I$ onto $D$ as follows: for $p\in \hat D_1\setminus I$ consider a point $p'\in\Pi_1^{-1}(p)$ and let ${\bf f}(p):=(\Pi\circ F\circ\theta^{-1}\circ\varphi\circ\theta)(p')$. By the construction of $F$, this definition is independent of the choice of $p'$. It is straightforward to prove that ${\bf f}$ is a biholomorphic map between $\overline{\hat D_1}\setminus I$ and $\overline{D}$. The domain $D$ is Kobayashi-hyperbolic as a biholomorphic image of the bounded domain $\hat D_1\setminus I$. 

It is shown in \cite{Kr} that a biholomorphic map between two hyperbolic Reinhardt domains in $\CC^n$ can be represented as the composition of their automorphisms and an elementary biholomorphic map between them. Since $\hat D_1\setminus I$ and $D$ do not intersect $I$, it follows from \cite{Kr} that all automorphisms of these domain are elementary. Therefore, ${\bf f}$ is an elementary map.

Further, $F_2^{-1}={\bf G}\circ F$, where ${\bf G}$ is an affine transformation of the form (\ref{affine}) with
$$
A=\left(
\begin{array}{ll}
a_1 & b_1\\
a_2 & b_2
\end{array}
\right).
$$
Hence $V:=\Pi(F(\tilde L_2))$ is mapped onto $U_2$ by an elementary map ${\bf g}$ of the form 
$$
\begin{array}{lll}
z&\mapsto&\hbox{const}\, z^{a_1}w^{b_1},\\
w&\mapsto&\hbox{const}\, z^{a_2}w^{b_2}.
\end{array}
$$
It is straightforward to verify that $f={\bf g}\circ{\bf f}$ on ${\bf f}^{-1}(V)\setminus(C_f\cup J_f)$, and therefore $f$ is an elementary map.

Assume now that $\tilde L_1$ is hypersurface (3). Then each of $\hat G_1$ and $\hat G_2$ has exactly two fixed points: 0 and $p_{\infty}$. Therefore either $\varphi(0)=0$, $\varphi(p_{\infty})=p_{\infty}$, or $\varphi(0)=p_{\infty}$, $\varphi(p_{\infty})=0$. Hence $\varphi$ preserves $B^2\cap {\cal L}_z$ and thus can be lifted to a holomorphic automorphism of $\tilde T_1$, that is, there exists a map $\tilde\varphi\in\hbox{Aut}(\tilde T_1)$ such that $\theta\circ\tilde\varphi=\varphi\circ\theta$. The map $\tilde\varphi$ is also defined on $\tilde L_1$ and can be chosen to satisfy the condition $\tilde\varphi(p_1')=p_2'$ which yields $\tilde\varphi\circ G_1\circ\tilde\varphi^{-1}\subset G_2$. Hence, as before, we can construct an intermediate hyperbolic domain $D$, a biholomorphic map ${\bf f}$ from $\hat D_1\setminus I$ onto $D$ that, as before, turns out to be elementary, and an elementary map ${\bf g}$ from a portion of $\partial D$ into $U_2$ such that $f={\bf g}\circ{\bf f}$ on a portion of $U_1$. Thus, we again obtain that $f$ is elementary.

The proof of the proposition is complete.\qed
\smallskip\\

It now remains to consider the cases when $\tilde L_2\subset \tilde L_1$ and $\tilde L_1$ is either hypersurface (2) or hypersurface (4) of (\ref{fourtypes}). As in the proof of Proposition  \ref{nottwofour}, we denote the map $\theta_1=\theta_2$ by $\theta$.  

Let first $\tilde L_1$ be hypersurface (2). Then the only fixed point of each of $\hat G_1$ and $\hat G_2$ in $\overline{B^2}$ is $p_{\infty}$, and therefore $\varphi$ has the form (\ref{formphi}). Assume that the group $\hat G_1$ contains an element $\hat g_1$ changing the $z$-coordinate. Then the only complex line preserved by $\hat g_1$ is ${\cal L}_z$ (see (\ref{group})). The map $\varphi\circ\hat g_1\circ\varphi^{-1}$ also preserves a unique complex line and it follows from (\ref{group}) that this line is also ${\cal L}_z$. Therefore, $\varphi$ preserves ${\cal L}_z$ (that is, we have $\zeta=0$). Arguing as in the last paragraph of the proof of Proposition \ref{nottwofour}, we obtain that $f$ in this case is elementary.

Assume now that none of the elements of $\hat G_1$ changes the $z$-coordinate, i.e., $\hat G_1$ consists of transformations of the form
\begin{equation}
\begin{array}{lll}
z&\mapsto& \displaystyle z,\\
\displaystyle w&\mapsto& w+i\alpha_1 n+i\beta_1 m, \qquad n,m\in\ZZ,
\end{array}\label{g11form}
\end{equation}
for some $\alpha_1,\beta_1\ge 0$, $\alpha_1+\beta_1> 0$. If in formula (\ref{formphi}) we have $\zeta=0$, then $\varphi$ preserves ${\cal L}_z$ and we again obtain that $f$ is elementary. Therefore, we shall assume that $\zeta\ne 0$. 

We shall show first of all that the group $\hat G_1$ has only one generator.

\begin{proposition}\label{onegenerator}\sl The group $\hat G_1$ consists of transformations of the form
\begin{equation}
\begin{array}{lll}
z&\mapsto& \displaystyle z,\\
\displaystyle w&\mapsto& w+i\alpha_0 n, \qquad n\in\ZZ,
\end{array}\label{formonegen}
\end{equation}
for some $\alpha_0>0$.
\end{proposition} 

\noindent {\bf Proof:} Let $\gamma:={\cal L}_z\cap S^3$ and $\gamma':=\varphi^{-1}(\gamma)$. Clearly, $\gamma'=\{z=-1/\lambda e^{-it}\zeta\}\cap S^3$. Let $\gamma_k':=\Bigl\{z=\ln_0(-1/\lambda e^{-it}\zeta)+i2\pi  k\Bigr\}\cap\tilde L_1$, for $k\in\ZZ$, be the curves in $\tilde L_1$ forming the set $\theta^{-1}(\gamma')$ (here $\ln_0$ denotes the principal branch of the logarithm). For some $k_0\in\ZZ$, $c\in\RR$ and sufficiently small $\varepsilon>0$ the circle $\tilde\gamma:=\Bigl\{\Bigl|z-\Bigl(\ln_0(-1/\lambda e^{-it}\zeta)+i2\pi k_0\Bigr)\Bigr|=\varepsilon\Bigr\}\cap\tilde L_1\cap\{\hbox{Im}\,w=c\}$ lies in $\tilde L_1\setminus\Bigl(\Pi_1^{-1}(C_f\cup J_f)\cup_{k\in\ZZ}\gamma_k'\Bigr)$. Recall that near $p_1'\in\tilde L_1\setminus\Pi_1^{-1}(C_f\cup J_f)$ we have
$$
\Pi_2\circ\eta\circ\varphi\circ\theta=f\circ\Pi_1,
$$ 
where $\eta$ is some analytic element of $\theta^{-1}$. The map in the right-hand side is well-defined everywhere on $\tilde L_1$, therefore, the analytic continuation of the map in the left-hand side along $\tilde\gamma$ produces a single-valued map. Clearly, after the analytic continuation of $\eta\circ\varphi\circ\theta$ along $\tilde\gamma$, its value changes by $(\pm 2\pi,0)$. Hence $G_2$ contains the map $\Lambda^z$ (defined at the beginning of this section).
Transformations in $G_2$ have the form
\begin{equation}
\left(
\begin{array}{l}
z\\
w
\end{array}
\right)
\mapsto
\left(
\begin{array}{l}
z\\
w
\end{array}
\right)
+i\left(
\begin{array}{l}
\alpha_2'\\
\alpha_2
\end{array}
\right)n+
i\left(
\begin{array}{l}
\beta_2'\\
\beta_2
\end{array}
\right)m,\quad n,m\in\ZZ,\label{g2form}
\end{equation}
for some linearly independent vectors $(\alpha_2',\alpha_2)$, $(\beta_2',\beta_2)\in\RR^2$. Since the map $\Lambda^z$ is contained in $G_2$, for some $n_0,m_0\in\ZZ$, it follows that
$$
\begin{array}{l}
\alpha_2'n_0+\beta_2'm_0=2\pi,\\  
\alpha_2n_0+\beta_2m_0=0.
\end{array}
$$
Hence $\alpha_2$ and $\beta_2$ are rationally dependent. 

Next, a straightforward calculation shows that the subgroup $\varphi\circ\hat G_1\circ\varphi^{-1}$ of $\hat G_2$ consists of the maps
\begin{equation}
\begin{array}{lll}
z&\mapsto& \displaystyle z,\\
\displaystyle w&\mapsto& w+i\lambda^2\alpha_1 n+i\lambda^2\beta_1 m, \qquad n,m\in\ZZ.
\end{array}\label{formg1}
\end{equation}
Taking into account that the general form of an element of $\hat G_2$ is
$$
\begin{array}{lll}
z&\mapsto& \exp(i\alpha_2'n+i\beta_2'm)z,\\
w&\mapsto& w+i\alpha_2n+i\beta_2m, \qquad n,m\in\ZZ,
\end{array}
$$
we see that $\alpha_1$ and $\beta_1$ are rationally dependent, and therefore transformations from $\hat G_1$ have the form (\ref{formonegen}) for some $\alpha_0>0$, as required.\qed
\smallskip\\

Let $D:=\left\{(z,w)\in\CC^2: |w|>\exp\left(|z|^2\right)\right\}$. We shall now construct a locally biholomorphic map ${\bf h}$ from $\hat D_1\setminus I$ onto $D\setminus I=D\setminus{\cal L}_z$.  
Obviously, the tube domain over the logarithmic diagram of $D$ is precisely $\tilde T_1$, and thus $\tilde T_1$ covers $D\setminus I$ by means of the map $\Pi$. We now use the map $\theta$ to construct a subgroup $\hat G\subset\hbox{Aut}(B^2)$ from the group $G$ acting on $\tilde T_1$ similarly to the way the groups $\hat G_1$, $\hat G_2$ were derived from $G_1$, $G_2$. Clearly, $\hat G$ consists of the transformations
$$
\begin{array}{lll}
z&\mapsto& z,\\
w&\mapsto& w+i2\pi n,\quad n\in\ZZ.
\end{array}
$$ 
Consider the following automorphism of $B^2$
$$
\varphi_1:\,\,
z\mapsto\delta_0 z,\quad w\mapsto\delta_0^2 w,
$$
where $\delta_0:=\sqrt{2\pi/\alpha_0}$. From (\ref{formonegen}) we obtain $\varphi_1\circ\hat G_1\circ\varphi_1^{-1}\subset\hat G$. For $p\in \hat D_1\setminus I$ consider a point $p'\in\Pi_1^{-1}(p)$, let $q\in \theta^{-1}\Bigl((\varphi_1\circ\theta)(p')\Bigr)$, and set ${\bf h}(p):=\Pi(q)$. Clearly, this definition is independent of the choices of $p'$ and $q$, and the map so defined is locally biholomorphic. Further, since $\varphi_1$ preserves ${\cal L}_z$, arguing as in the last paragraph of the proof of Proposition \ref{nottwofour}, we can show that ${\bf h}$ is elementary.

We shall now pause to describe the general form of a bounded domain whose complement to $I$ can be mapped onto $D\setminus I$ by means of an elementary map, as well as such elementary maps. For $a_1,b_1,c_1,d_1\in\ZZ$, $a_1>0$, $b_1>0$, $c_1\ge 0$, $d_1> 0$, such that $a_1d_1-b_1c_1>0$, and $C_1>0$, $E_1>0$ consider the domain
$$
\begin{array}{l}
R(a_1,b_1,c_1,d_1,C_1,E_1):=\\
\hspace{2cm}\left\{(z,w)\in\CC^2: C_1|z|^{c_1}|w|^{-d_1}>\exp\left(E_1 |z|^{2a_1}|w|^{-2b_1}\right),\,w\ne 0\right\}.
\end{array}
$$
The general form of an elementary map from $R(a_1,b_1,c_1,d_1,C_1,E_1)\setminus I$ onto $D\setminus I$ is
\begin{equation}
\begin{array}{lll}
z&\mapsto& e^{i\tau_1}\sqrt{C_1} z^{a_1}w^{-b_1},\\
\vspace{0mm}&&\\
w&\mapsto& e^{i\tau_2}\sqrt{E_1} z^{c_1}w^{-d_1},
\end{array}\label{algmap2}
\end{equation}
where $\tau_1,\tau_2\in\RR$. We observe that $R(a_1,b_1,c_1,d_1,C_1,E_1)\cap {\cal L}_z\ne\emptyset$ only if $c_1=0$. It is straightforward to show that a bounded domain whose complements to $I$ can be mapped onto $D\setminus I$ by an elementary map, up to permutation of the variables, is some $R(a_1,b_1,c_1,d_1,C_1,E_1)$ minus a closed subset of ${\cal L}_z$. Since $\hat D_1$ is pseudoconvex, up to permutation of the variables, we have either $\hat D_1=R(a_1,b_1,c_1,d_1,C_1,C_2)$ or $\hat D_1=R(a_1,b_1,0,d_1,C_1,C_2)\setminus {\cal L}_z$, for some $a_1,b_1,c_1,d_1, C_1, E_1$.

We shall now construct a biholomorphic map ${\bf f}:D\ra D^{\lambda}$, where $D^{\lambda}:=\left\{(z,w)\in\CC^2: |w|^{\lambda^2}>\exp\left(|z|^2\right)\right\}$. The tube domain
$$
T^{\lambda}:=\left\{(z,w)\in\CC^2:\lambda^2\hbox{Re}\,w>\exp\left(2\hbox{Re}\,z\right)\right\}
$$ 
covers $D^{\lambda}\setminus I$ by means of $\Pi$ and is mapped into $B^2$ by the map
$$
\theta^{\lambda}:\,\,
z\mapsto e^z,\quad w\mapsto\lambda^2 w.
$$
Denote by $\hat G^{\lambda}\subset\hbox{Aut}(B^2)$ the subgroup obtained by means of $\theta^{\lambda}$ from the group $G$ acting on $T^{\lambda}$. Clearly, $\hat G^{\lambda}$ consists of the transformations
$$
\begin{array}{lll}
z&\mapsto& z,\\
w&\mapsto& w+i\lambda^2 2\pi n,\quad n\in\ZZ.
\end{array}
$$ 

Let $\varphi_2$ be the following automorphism of $B^2$
$$
\begin{array}{lll}
z&\mapsto&\lambda e^{it}z+\delta_0\zeta,\\
w&\mapsto&\lambda^2w+2\delta_0\overline{\zeta}\lambda e^{it}z+\delta_0^2|\zeta|^2+i\delta_0\mu.
\end{array}
$$
A straightforward calculation shows that $\varphi_2\circ\hat G\circ\varphi_2^{-1}=\hat G^{\lambda}$. Let ${\cal L}':=\left\{z=-1/\lambda e^{-it}\delta_0\zeta\right\}$. For $p\in D\setminus ({\cal L}_z\cup{\cal L}')$ consider a point $p'\in\Pi^{-1}(p)$, let $q\in\theta^{{\lambda}^{-1}}\Bigl((\varphi_2\circ \theta)(p')\Bigr)$, and set ${\bf f}(p):=\Pi(q)$. This definition is clearly independent of the choices of $p'$ and $q$. It is straightforward to verify that ${\bf f}$ maps $D\setminus ({\cal L}_z\cup{\cal L}')$ biholomorphically onto $D^{\lambda}\setminus ({\cal L}_z\cup{\cal L}'')$, where ${\cal L}'':=\{z=\delta_0\zeta\}$. Since $D$ and $D^{\lambda}$ have bounded realizations, ${\bf f}$ extends to a map (also denoted by ${\bf f}$) from $D$ onto $D^{\lambda}$. This map is biholomorphic and ${\bf f}(D\cap{\cal L}_z)=D^{\lambda}\cap{\cal L}''$, ${\bf f}(D\cap{\cal L}')=D^{\lambda}\cap{\cal L}_z$. Further, ${\bf f}$ can be represented as ${\bf f}={\bf f}_1\circ{\bf f}_2$ with
$$
{\bf f}_1:\,\, 
z\mapsto \lambda z,\quad w\mapsto w,
$$
and ${\bf f}_2\in\hbox{Aut}(D)$. The map ${\bf f}_2$ has the form (see \cite{Kr}, \cite{Sh})
\begin{equation}
\begin{array}{lll}
z&\mapsto& e^{i\tau_1}z+s,\\
w&\mapsto&
e^{i\tau_2}\exp\left(2\overline{s}e^{i\tau_1}z+
|s|^2\right)w,
\end{array}\label{formmf}
\end{equation}
where $\tau_1,\tau_2\in\RR$, $s\in\CC^*$.

Finally, we define a locally biholomorphic map ${\bf g}$ from $D^{\lambda}\setminus I=D^{\lambda}\setminus{\cal L}_z$ onto $\Omega:=\Pi_2(\tilde T_1)$. It is constructed similarly to the map ${\bf h}$. Consider the following automorphism of $B^2$
$$
\varphi_3:\,\,
z\mapsto\frac{1}{\delta_0}z,\quad w\mapsto\frac{1}{\delta^2_0}w.
$$
It follows from (\ref{formonegen}), (\ref{formg1}) that $\varphi_3\circ\hat G^{\lambda}\circ\varphi_3^{-1}\subset\hat G_2$. For $p\in D^{\lambda}\setminus I$ let $p'\in\Pi^{-1}(p)$, $q\in \theta^{-1}\Bigl((\varphi_3\circ\theta^{\lambda})(p')\Bigr)$, and set ${\bf g}(p):=\Pi_2(q)$. As before, this definition is independent of the choices of $p'$ and $q$ (recall that $G_2$ contains the transformation $\Lambda^z$), and the map so defined is locally biholomorphic. Since $\varphi_3$ preserves ${\cal L}_z$, arguing again as in the last paragraph of the proof of Proposition \ref{nottwofour}, we obtain that ${\bf g}$ is elementary.

The composition ${\bf g}\circ{\bf f}\circ{\bf h}$ maps $V:={\bf h}^{-1}\Bigl(({\bf g}\circ{\bf f})^{-1}(U_2)\setminus I\Bigr) $ into $U_2\subset\partial\Omega$. Since $\varphi=\varphi_3\circ\varphi_2\circ\varphi_1$, it follows that $f={\bf g}\circ{\bf f}\circ{\bf h}$ on $V\setminus(C_f\cup J_f)$. Therefore, $f={\bf g}\circ{\bf f}\circ{\bf h}$ on $\hat D_1\setminus {\bf h}^{-1}({\cal L}')$. Clearly, $f$ maps $\hat D_1\setminus {\bf h}^{-1}({\cal L}')$ onto a set of the form $\Omega\setminus U$, where either $U=\emptyset$ (if $\hat D_1\cap I\ne\emptyset$) or $U={\bf g}\left({\cal L}''\cap D^{\lambda}\right)$ (if $\hat D_1\cap I=\emptyset$). 

If $U\ne\emptyset$, then $f(\hat D_1)$ is not a Reinhardt domain, because $s\ne 0$ in formula (\ref{formmf}). This shows that in fact $U=\emptyset$, that is, $\hat D_1\cap I\ne\emptyset$ which implies that, up to permutation of the variables, $\hat D_1=R(a_1,b_1,0,d_1,C_1,E_1)$ for some $a_1,b_1,c_1,d_1,C_1,E_1$, and ${\bf h}$ has the form (\ref{algmap2}).

Further, $\Omega$ is a bounded Reinhardt domain not intersecting $I$, and $D^{\lambda}\setminus I$ is mapped onto $\Omega$ by an elementary map. It is not difficult to describe all such domains and the corresponding elementary maps. A domain of this kind has the form
$$
\left\{(z,w)\in\CC^2: C_2|z|^{\frac{c_2}{\Delta}}|w|^{\frac{a_2}{\Delta}}>\exp\left(E_2 |z|^{\frac{2d_2}{\Delta}}|w|^{\frac{2b_2}{\Delta}}\right),
z\ne 0,\, w\ne 0\right\},
$$
where $\Delta:=a_2d_2-b_2c_2$, $\Delta\ne 0$, $a_2,b_2,c_2,d_2\in\ZZ$, $a_2\ge 0$, $b_2>0$, $c_2\le 0$, $d_2<0$, $C_2>0$, $E_2>0$. The general form of an elementary map from $D^{\lambda}\setminus I$ onto the above domain is
\begin{equation}
\begin{array}{lll}
z&\mapsto& \hbox{const}\, z^{a_2} w^{-b_2},\\
w&\mapsto& \hbox{const}\, z^{-c_2} w^{d_2}.
\end{array}\label{formmapg}
\end{equation}
In particular, $\Omega$ and ${\bf g}$ must have these forms.

Since $U=\emptyset$, we obtain $\hat D_2=\Omega\cup{\bf g}\left({\cal L}_z\cap D^{\lambda}\right)$. If $a_2>0$ and $c_2<0$, it follows from (\ref{formmapg}) that ${\bf g}\left({\cal L}_z\cap D^{\lambda}\right)=\{0\}$. However, $\Omega\cup\{0\}$ is not an open set in this case, and therefore either $a_2=0$ or $c_2=0$. If $a_2=0$, then $c_2<0$ and we have
$$
\hat D_2=\left\{(z,w)\in\CC^2: C_2|z|<\exp\left(-E_2' |z|^{-\frac{2d_2}{b_2c_2}}|w|^{-\frac{2}{c_2}}\right),\,z\ne 0\right\}    
$$
for some $E_2'>0$; if $c_2=0$, then $a_2>0$ and we have
$$
\hat D_2=\left\{(z,w)\in\CC^2: C_2|w|<\exp\left(-E_2'' |z|^{\frac{2}{a_2}}|w|^{\frac{2b_2}{a_2d_2}}\right),\,w\ne 0\right\},    
$$
for some $E_2''>0$. The above two classes of domains are obtained from one another by permutation of the variables.

It is clear that every subdomain of $\hat D_1$ mapped properly by $f$ onto a bounded Reinhardt domain and whose envelope of holomorphy coincides with $\hat D_1$ up to permutation of the variables has the form
$$
\begin{array}{ll}
\displaystyle \Biggl\{(z,w)\in\CC^2: C_1^{*\frac{1}{d_1}}\exp\left(-\frac{E_1}{d_1} |z|^{\frac{2}{a_1}}|w|^{-2b_1}\right)&<|w|<\\
&\displaystyle C_1^{\frac{1}{d_1}}\exp\left(-\frac{E_1}{d_1} |z|^{\frac{2}{a_1}}|w|^{-2b_1}\right)\Biggr\},
\end{array}
$$
for some $0\le C_1^*<C_1$, and hence $D_1$ is of this form. We thus have obtained (iv) of Theorem \ref{mainresult}.

Let now $\tilde L_1$ be hypersurface (4) of (\ref{fourtypes}). For the purposes of this case we realize $S^3$ as $\{(z,w)\in\CC^2: |z|^2+|w|^2=1\}$ and $B^2$ as $\{(z,w)\in\CC^2: |z|^2+$\linebreak $|w|^2<1\}$. Then we have $\theta=\Pi$, $\theta(\tilde T_1)=B^2\setminus I$, and each of $\hat G_1$, $\hat G_2$ consists of transformations of the form
$$
\begin{array}{lll}
z&\mapsto& e^{i\alpha_1}z,\\
w&\mapsto& e^{i\alpha_2}w,
\end{array}
$$
where the vector $(\alpha_1,\alpha_2)$ varies over a lattice in $\RR^2$.

Assume first that $\varphi(B^2\cap I)=B^2\cap I$. In this case $\varphi$ can be lifted to an automorphism $\tilde\varphi$ of $\tilde T_1$ such that $\tilde\varphi\circ G_1\circ\tilde\varphi^{-1}\subset G_2$. Then, arguing as in the last paragraph of the proof of Proposition \ref{nottwofour}, we see that $f$ is elementary.

Assume now that $\varphi(B^2\cap I)\ne B^2\cap I$ and $\varphi(B^2\cap{\cal L}_w)=B^2\cap{\cal L}_w$. Then $\varphi$ has the form
$$
\begin{array}{lll}
z&\mapsto&\displaystyle e^{it_1}\frac{z-a}{1-\overline{a}z},\\
\vspace{0mm}&&\\
w&\mapsto&\displaystyle e^{it_2}\frac{\sqrt{1-|a|^2}}{1-\overline{a}z}w,
\end{array}
$$
where $|a|<1$, $a\ne 0$, and $t_1,t_2\in\RR$. It is now clear from the inclusion $\varphi\circ\hat G_1\circ\varphi^{-1}\subset\hat G_2$ that none of the elements of $\hat G_1$ changes the $z$-coordinate.

Consider the tube domain
$$
T:=\left\{(z,w)\in\CC^2:\hbox{Re}\,z>\exp\left(2\hbox{Re}\,w\right)\right\}.
$$
This domain covers $B^2\setminus{\cal L}_w$ by means of the map 
$$
\hat\theta:\,\,
z\mapsto\frac{z-1}{z+1},\quad w\mapsto -\frac{2e^w}{z+1},
$$
and $\tilde T_1$ covers $T\setminus\{z=1\}$ by means of the map 
$$
\check{\theta}:\,\,
z\mapsto -\frac{e^z+1}{e^z-1},\quad w\mapsto w-\ln_0(e^z-1).
$$
Clearly, $\theta=\hat\theta\circ\check{\theta}$. Since the groups $\hat G_1$, $\hat G_2$ preserve ${\cal L}_w$, their elements can be lifted to automorphisms of $T$. Similarly, $\varphi$ can be lifted to an automorphism of $T$. The general form of a lift of $\varphi$ is
\begin{equation}
\begin{array}{lll}
z&\mapsto&\displaystyle-\frac{(e^{it_1}(1-a)+1-\overline{a})z-e^{it_1}(1+a)+1+\overline{a}}{(e^{it_1}(1-a)-1+\overline{a})z-e^{it_1}(1+a)-1-\overline{a}},\\
\vspace{0mm}&&\\
w&\mapsto&\displaystyle w+\ln\frac{-2e^{it_2}\sqrt{1-|a|^2}}{(e^{it_1}(1-a)-1+\overline{a})z-e^{it_1}(1+a)-1-\overline{a}},
\end{array}\label{liftphi}
\end{equation}
where $\ln$ is a branch of the logarithm.

For arbitrary $g_1\in G_1$ consider the locally defined self-map $\tilde g_1=\check{\theta}\circ g_1\circ \check{\theta}^{-1}$ of $T$. Clearly, it coincides with a lift of $\hat g_1\in \hat G_1$, and hence extends to an automorphism of $T$. Let $\tilde G_1:=\{\tilde g_1, g_1\in G_1\}$, and let $\tilde G_2$ be the group constructed from $G_2$ in the same way. The groups $\tilde G_1$, $\tilde G_2$ are Abelian and have at most two generators. Observe that a lift $\tilde\varphi$ of $\varphi$ to an automorphism of $T$ can be chosen so that $\tilde\varphi\circ\tilde G_1\circ\tilde\varphi^{-1}\subset\tilde G_2$.

The group $G_1$ consist of transformations of the form
\begin{equation}
\left(
\begin{array}{l}
z\\
w
\end{array}
\right)
\mapsto
\left(
\begin{array}{l}
z\\
w
\end{array}
\right)
+i\left(
\begin{array}{l}
\alpha_1'\\
\alpha_1
\end{array}
\right)n+
i\left(
\begin{array}{l}
\beta_1'\\
\beta_1
\end{array}
\right)m,\quad n,m\in\ZZ,\label{groupg1forrrrm}
\end{equation}
for some linearly independent vectors $(\alpha_1',\alpha_1), (\beta_1',\beta_1)\in\RR^2$. Then $\hat G_1$ consists of the maps
$$
\begin{array}{lll}
z&\mapsto&\exp\left(i\alpha_1'n+i\beta_1'm\right)z,\\   
w&\mapsto&\exp\left(i\alpha_1n+i\beta_1m\right)w, \quad n,m\in\ZZ.
\end{array}
$$
Since no map in $\hat G_1$ changes the $z$-coordinate, it follows that $\alpha_1',\beta_1'\in 2\pi\cdot\ZZ$. Therefore, elements of $\tilde G_1$ have the form (\ref{g11form}). It then follows from (\ref{liftphi}) that every element of $\tilde G_1$ commutes with every lift of $\varphi$ to an automorphism of $T$. Hence $\tilde G_1\subset \tilde G_2$.

Next, if the group $G_2$ is given by (\ref{g2form}), arguing as in the proof of Proposition \ref{onegenerator}, we obtain that $G_2$ contains the map $\Lambda^z$. As before, this yields that $\alpha_2$ and $\beta_2$ are rationally dependent. Further, transformations from $\tilde G_2$ have the following form
$$
\begin{array}{lll}
z &\mapsto& \displaystyle\frac{(1+C(n,m))z+1-C(n,m)}{(1-C(n,m))z+1+C(n,m)},\\
\vspace{0mm}&&\\
w&\mapsto & \displaystyle w+\ln\frac{2}{(1-C(n,m))z+1+C(n,m)}+i\alpha_2n+i\beta_2m,\quad n,m\in\ZZ,
\end{array}
$$
where $C(n,m):=\exp\left(i\alpha_2'n+i\beta_2'm\right)$. For elements of $\tilde G_1$ the corresponding constants $C(n,m)$ are necessarily equal to 1, which implies that $\alpha_1$ and $\beta_1$ are also rationally dependent. Therefore, transformations from $\tilde G_1$ have the form (\ref{formonegen}) for some $\alpha_0>0$.

Let $D^{\alpha_0}:=\left\{(z,w)\in\CC^2: |z|^2+|w|^{\frac{\alpha_0}{\pi}}<1\right\}$. We shall now construct a locally biholomorphic map ${\bf h}$ from $\hat D_1\setminus I$ onto $D^{\alpha_0}\setminus I $. Clearly, $\tilde T_1$ covers $D^{\alpha_0}\setminus I$ by means of the map
$$
\Pi^{\alpha_0}:\,\, 
z\mapsto e^z,\quad w\mapsto e^{\frac{2\pi}{\alpha_0}w}.
$$
The group $G^{\alpha_0}$ constructed from $D^{\alpha_0}$ in the same way as $G_1$ and $G_2$ were constructed from $\hat D_1$ and $\hat D_2$, consists of the following transformations
$$
\begin{array}{lll}
z&\mapsto& z+i2\pi n,\\
w&\mapsto& w+i\alpha_0 m,\qquad n,m\in\ZZ.
\end{array}
$$
For $p\in \hat D_1\setminus I$ consider a point $p'\in\Pi_1^{-1}(p)$ and set ${\bf h}(p):=\Pi^{\alpha_0}(p')$. Clearly, this definition is independent of the choice of $p'$, and the map so defined is locally biholomorphic. The automorphism of $B^2$ that corresponds to ${\bf h}$ is the identity, and therefore preserves $I$. Arguing as in the last paragraph of the proof of Proposition \ref{nottwofour}, we see that ${\bf h}$ is elementary.

It is straightforward to describe the general form of a bounded domain not intersecting $I$ that can be mapped onto $D^{\alpha_0}\setminus I$ by an elementary map, as well as such elementary maps. Such a domain must have the following form
\begin{equation}
\left\{(z,w)\in\CC^2: C_1|z|^{2a_1}|w|^{2b_1}+E_1 |z|^{\frac{\alpha_0 c_1}{\pi}}|w|^{\frac{\alpha_0 d_1}{\pi}}<1,\,z\ne 0,\, w\ne 0\right\},\label{domainv}
\end{equation}
where $a_1,b_1,c_1,d_1\in\ZZ$, with either $a_1d_1-b_1c_1>0$, $a_1\ge 0$, $b_1\le 0$, $c_1\le 0$, $d_1\ge 0$, or $a_1d_1-b_1c_1<0$, $a_1\le 0$, $b_1\ge 0$, $c_1\ge 0$, $d_1\le 0$, and $C_1>0$, $E_1>0$. An elementary map that takes domain (\ref{domainv}) onto $D^{\alpha_0}\setminus I$ has the form
$$
\begin{array}{lll}
z&\mapsto& e^{i\tau_1}\sqrt{C_1}z^{a_1}w^{b_1},\\
\vspace{0mm}&&\\
w&\mapsto& \displaystyle e^{i\tau_2}E_1^{\frac{\pi}{\alpha_0}}z^{c_1}w^{d_1},
\end{array}
$$
where $\tau_1,\tau_2\in\RR$. Thus, $\hat D_1\setminus I$ and ${\bf h}$ must have the forms described above. Since $\hat D_1$ is pseudoconvex, it is either domain (\ref{domainv}) or, up to permutation of the variables, one of the following domains
$$
\begin{array}{l}
\hspace{-1cm}\left\{(z,w)\in\CC^2: C_1|w|^{2b_1}+E_1 |z|^{\frac{\alpha_0 c_1}{\pi}}|w|^{\frac{\alpha_0 d_1}{\pi}}<1,\, w\ne 0\right\},\\
\hspace{-1cm}\hbox{(here $a_1=0$, $b_1>0$, $c_1>0$, $d_1\le 0$)},
\end{array}
$$
\begin{equation}
\begin{array}{l}
\left\{(z,w)\in\CC^2: C_1|z|^{2a_1}|w|^{2b_1}+E_1 |w|^{\frac{\alpha_0 d_1}{\pi}}<1,\, w\ne 0\right\},\\
\hbox{(here $a_1>0$, $b_1\le 0$, $c_1=0$, $d_1>0$)},
\end{array}\label{formsd_12}
\end{equation}
\begin{equation}
\begin{array}{l}
\hspace{-3.1cm}\left\{(z,w)\in\CC^2: C_1|z|^{2a_1}+E_1 |w|^{\frac{\alpha_0 d_1}{\pi}}<1,\right\},\\
\hspace{-3.1cm}\hbox{(here $a_1>0$, $b_1=0$, $c_1=0$, $d_1>0$)},
\end{array}\label{formsd_1}
\end{equation}
for some $a_1,b_1,c_1,d_1,C_1,E_1$.

We shall now construct ${\bf f}\in\hbox{Aut}(D^{\alpha_0})$. Let ${\cal L}':=\{z=a\}$ and ${\cal L}'':=\{z=-e^{it_1}a\}$. It is straightforward to observe that $\tilde G^{\alpha_0}=\tilde G_1$. In particular, elements of $\tilde G^{\alpha_0}$ commute with $\tilde\varphi$, which yields $\tilde\varphi\circ\tilde G^{\alpha_0}\circ\tilde\varphi^{-1}=\tilde G^{\alpha_0}$. For $p\in D^{\alpha_0}\setminus (I\cup{\cal L}')$ consider a point $p'\in\Pi^{{\alpha_0}^{-1}}(p)$, let $q\in\check{\theta}^{-1}\Bigl((\tilde\varphi\circ \check{\theta})(p')\Bigr)$, and set ${\bf f}(p):=\Pi^{\alpha_0}(q)$. Clearly, this definition is independent of the choices of $p'$ and $q$. It is straightforward to verify that ${\bf f}$ maps $D^{\alpha_0}\setminus (I\cup{\cal L}')$ biholomorphically onto $D^{\alpha_0}\setminus (I\cup {\cal L}'')$. Since $D^{\alpha_0}$ is bounded, ${\bf f}$ extends to a map (that we also denote by ${\bf f}$) from $D^{\alpha_0}$ onto itself. This map is biholomorphic and ${\bf f}(D^{\alpha_0}\cap{\cal L}')\subset D^{\alpha_0}\cap I$, ${\bf f}(D^{\alpha_0}\cap I)\subset D^{\alpha_0}\cap (I\cup {\cal L}'')$. It is now clear that ${\bf f}$ has the form
\begin{equation}
\begin{array}{lll}
z&\mapsto&\displaystyle e^{it_1}\frac{z-a}{1-\overline{a}z},\\
\vspace{0mm}&&\\
w&\mapsto&\displaystyle e^{it}\frac{(1-|a|^2)^{\frac{\pi}{\alpha_0}}}{(1-\overline{a}z)^{\frac{2\pi}{\alpha_0}}}w,
\end{array}\label{formmmf}
\end{equation}
where $t\in\RR$.

Finally, we define a locally biholomorphic map ${\bf g}$ from $D^{\alpha_0}\setminus I$ onto $\Omega:=\Pi_2(\tilde T_1)$. It is constructed similarly to the map ${\bf h}$. For $p\in \hat D^{\alpha_0}\setminus I$ consider a point $p'\in\Pi^{\alpha_0^{-1}}(p)$ and set ${\bf g}(p):=\Pi_2(p')$. Since $\tilde G_1\subset \tilde G_2$ and the map $\Lambda^z$ belongs to $G_2$, the map 
$$
\begin{array}{lll}
z&\mapsto& z,\\
w&\mapsto& w+i\alpha_0,
\end{array}
$$
belongs to $G_2$ as well. Therefore, the above definition of ${\bf g}$ is independent of the choice of $p'$. The map so defined is locally biholomorphic.  The automorphism of $B^2$ corresponding to ${\bf g}$ is the identity, and therefore preserves $I$. Arguing as in the last paragraph of the proof of Proposition \ref{nottwofour}, we see that ${\bf g}$ is elementary.

The composition ${\bf g}\circ{\bf f}\circ{\bf h}$ maps $V:={\bf h}^{-1}\Bigl(({\bf g}\circ{\bf f})^{-1}(U_2)\setminus I\Bigr) $ into $U_2\subset\partial\Omega$. It is straightforward to verify that $f={\bf g}\circ{\bf f}\circ{\bf h}$ on $V\setminus(C_f\cup J_f)$. Therefore, $f={\bf g}\circ{\bf f}\circ{\bf h}$ on $\hat D_1\setminus {\bf h}^{-1}({\cal L}'\cup{\cal L}_w)$. Clearly, $f$ maps $\hat D_1\setminus {\bf h}^{-1}({\cal L}'\cup{\cal L}_w)$ onto $\Omega\setminus U$, where either $U=\emptyset$, or $U={\bf g}(D^{\alpha_0}\cap{\cal L}'')$. 

If $U\ne\emptyset$, then $f(\hat D_1)$ is not a Reinhardt domain,
since $a\ne 0$ in formula (\ref{formmmf}). Hence in fact $U=\emptyset$, that is, ${\bf h}(\hat D_1)\cap{\cal L}_z$ contains the punctured disc $\{z=0,\,0<|w|<1\}$, and therefore, up to permutation of the variables, $\hat D_1$ has one of the forms (\ref{formsd_12}), (\ref{formsd_1}).

Further, $\Omega$ is a bounded Reinhardt domain not intersecting $I$, and $D^{\alpha_0}\setminus I$ can be mapped onto $\Omega$ by an elementary map. It is not hard to describe the general form of such domains and elementary maps. A domain of this kind has the form
\begin{equation}
\begin{array}{l}
\Bigl\{(z,w)\in\CC^2: C_2|z|^{\frac{2d_2}{\Delta}}|w|^{-\frac{2b_2}{\Delta}}+E_2|z|^{-\frac{\alpha_0c_2}{\pi\Delta}}|w|^{\frac{\alpha_0 a_2}{\pi\Delta}}<1,\\ \hspace{8cm}z\ne 0,\, w\ne 0\Bigr\},
\end{array}\label{onemoreform}
\end{equation}
where $\Delta:=a_2d_2-b_2c_2$, $\Delta\ne 0$, $a_2,b_2,c_2,d_2\in\ZZ$, $a_2\ge 0$, $b_2\ge0$, $c_2\ge 0$, $d_2\ge 0$, $C_2>0$, $E_2>0$. An elementary map taking $D^{\alpha_0}\setminus I$ onto this domain is of the form
\begin{equation}
\begin{array}{lll}
z&\mapsto& \hbox{const}\, z^{a_2} w^{b_2},\\
w&\mapsto& \hbox{const}\, z^{c_2} w^{d_2}.
\end{array}\label{formmmmmmg}
\end{equation}
In particular, $\Omega$ and ${\bf g}$ must have these forms.

Assume that $\hat D_1$ has the form (\ref{formsd_1}). Since $U=\emptyset$, we obtain $\hat D_2=\Omega\cup{\bf g}\left(D^{\alpha_0}\cap I\right)$.
If either $a_2>0$ and $c_2>0$, or $b_2>0$ and $d_2>0$, it follows from (\ref{formmmmmmg}) that ${\bf g}\left(D^{\alpha_0}\cap I\right)=\{0\}$. However, $\Omega\cup\{0\}$ is not an open set in this case, and hence either $a_2=0$, $d_2=0$, or $b_2=0$, $c_2=0$. In the first case $b_2>0$, $c_2>0$ and
\begin{equation}
\hat D_2=\left\{(z,w)\in\CC^2: C_2|w|^{\frac{2}{c_2}}+E_2 |z|^{\frac{\alpha_0}{\pi b_2}}<1\right\}.\label{domm}    
\end{equation}
We then have either $D_1=\hat D_1$, $D_2=\hat D_2$, or 
\begin{equation}
\begin{array}{ll}
\displaystyle D_1=\Biggl\{(z,w)\in\CC^2:& C_1|z|^{2a_1}<1,
E_1^{*-\frac{\pi}{\alpha_0 d_1}}\left(1-C_1|z|^{2a_1}\right)^{\frac{\pi}{\alpha_0 d_1}}<\\
&|w|<\displaystyle E_1^{-\frac{\pi}{\alpha_0 d_1}}\left(1-C_1|z|^{2a_1}\right)^{\frac{\pi}{\alpha_0 d_1}}\Biggr\},
\end{array}\label{unexp1}
\end{equation}
$$
\begin{array}{ll}
\displaystyle D_2=\Biggl\{(z,w)\in\CC^2:& C_2|w|^{\frac{2}{c_2}}<1,
E_2^{*-\frac{\pi b_2}{\alpha_0}}\left(1-C_2|w|^{\frac{2}{c_2}}\right)^{\frac{\pi b_2}{\alpha_0}}<\\
&\displaystyle |z|<E_2^{-\frac{\pi b_2}{\alpha_0}}\left(1-C_2|w|^{\frac{2}{c_2}}\right)^{\frac{\pi b_2}{\alpha_0}}\Biggr\},
\end{array}
$$
for some $E_1<E_1^*\le\infty$, $E_2<E_2^*\le\infty$. 
Similarly, in the second case $a_2>0$, $d_2>0$, and
\begin{equation}
\hat D_2=\left\{(z,w)\in\CC^2: C_2|z|^{\frac{2}{a_2}}+E_2 |w|^{\frac{\alpha_0}{\pi d_2}}<1\right\}.\label{dommmmmm} 
\end{equation}
We then have either $D_1=\hat D_1$, $D_2=\hat D_2$, or $D_1$ has the form (\ref{unexp1}) and
$$
\begin{array}{ll}
\displaystyle D_2=\Biggl\{(z,w)\in\CC^2:& C_2|z|^{\frac{2}{a_2}}<1,
E_2^{*-\frac{\pi d_2}{\alpha_0}}\left(1-C_2|z|^{\frac{2}{a_2}}\right)^{\frac{\pi d_2}{\alpha_0}}<|w|<\\
&\displaystyle E_2^{-\frac{\pi d_2}{\alpha_0}}\left(1-C_2|z|^{\frac{2}{a_2}}\right)^{\frac{\pi d_2}{\alpha_0}}\Biggr\},
\end{array}
$$
for some $E_2<E_2^*\le\infty$. The above two forms of $\hat D_2$ are obtained from one another by permutation of the variables. 

Assume now that $\hat D_1$ has the form (\ref{formsd_12}). Then, as before, either $a_2=0$, $c_2>0$, or $a_2>0$, $c_2=0$. In the first case $b_2>0$ and
$$
\hat D_2=\left\{(z,w)\in\CC^2: C_2 |z|^{-\frac{2d_2}{b_2c_2}}|w|^{\frac{2}{c_2}}+E_2 |z|^{\frac{\alpha_0}{\pi b_2}}<1, z\ne 0\right\},
$$
in the second case $d_2>0$ and 
$$
\hat D_2=\left\{(z,w)\in\CC^2: C_2|z|^{\frac{2}{a_2}}|w|^{-\frac{2b_2}{a_2d_2}}+E_2 |w|^{\frac{\alpha_0}{\pi d_2}}<1, w\ne 0\right\}.
$$
The above two domains are obtained from one another by permutation of the variables. In the first case we obtain
\begin{equation}
\begin{array}{ll}
\displaystyle D_1=\Biggl\{(z,w)\in\CC^2: & \displaystyle C_1|z|^{2a_1}|w|^{2b_1}<1,\\
&\displaystyle
E_1^{*-\frac{\pi}{\alpha_0 d_1}}\left(1-C_1|z|^{2a_1}|w|^{2b_1}\right)^{\frac{\pi}{\alpha_0 d_1}}< |w|<\\
&\displaystyle E_1^{-\frac{\pi}{\alpha_0 d_1}}\left(1-C_1|z|^{2a_1}|w|^{2b_1}\right)^{\frac{\pi}{\alpha_0 d_1}}\Biggr\},
\end{array}\label{unexp2}
\end{equation}
$$
\begin{array}{ll}
\displaystyle D_2=\Biggl\{(z,w)\in\CC^2:& \displaystyle C_2|z|^{-\frac{2d_2}{b_2c_2}}|w|^{\frac{2}{c_2}}<1,\\
&\displaystyle E_2^{*-\frac{\pi b_2}{\alpha_0}}\left(1-C_2|z|^{-\frac{2d_2}{b_2c_2}}|w|^{\frac{2}{c_2}}\right)^{\frac{\pi b_2}{\alpha_0}}<|z|<\\
&\displaystyle E_2^{-\frac{\pi b_2}{\alpha_0}}\left(1-C_2|z|^{-\frac{2d_2}{b_2c_2}}|w|^{\frac{2}{c_2}}\right)^{\frac{\pi b_2}{\alpha_0}}\Biggr\},
\end{array}
$$
for some $E_1<E_1^*\le\infty$, $E_2<E_2^*\le\infty$. Similarly, in the second case $D_1$ has the form (\ref{unexp2}) and
$$
\begin{array}{ll}
\displaystyle D_2=\Biggl\{(z,w)\in\CC^2:& \displaystyle C_2|z|^{\frac{2}{a_2}}|w|^{-\frac{2b_2}{a_2d_2}}<1,\\
&\displaystyle E_2^{*-\frac{\pi d_2}{\alpha_0}}\left(1-C_2|z|^{\frac{2}{a_2}}|w|^{-\frac{2b_2}{a_2d_2}}\right)^{\frac{\pi d_2}{\alpha_0}}<|w|<\\
&\displaystyle E_2^{-\frac{\pi d_2}{\alpha_0}}\left(1-C_2|z|^{\frac{2}{a_2}}|w|^{-\frac{2b_2}{a_2d_2}}\right)^{\frac{\pi d_2}{\alpha_0}}\Biggr\},
\end{array}
$$
for some $E_2<E_2^*\le\infty$.

Similar considerations in the cases when $\varphi(B^2\cap I)\ne B^2\cap I$ and either $\varphi(B^2\cap{\cal L}_z)=B^2\cap{\cal L}_z$ or $\varphi(B^2\cap{\cal L}_z)=B^2\cap{\cal L}_w$ or $\varphi(B^2\cap{\cal L}_w)=B^2\cap{\cal L}_z$ lead to the same descriptions of $D_1$, $D_2$, and $f$. We thus have obtained (v) of Theorem \ref{mainresult}.

Assume finally that $\varphi(B^2\cap{\cal L}_z)\not\subset I$ and $\varphi(B^2\cap{\cal L}_w)\not\subset I$. Arguing as in the proof of Proposition \ref{onegenerator}, we can prove that in $G_2$ contains the map $\Lambda^z$ as well as the map $\Lambda^w$ (see the beginning of the section for definitions). Therefore, all elements of the inverse of the matrix $A$ corresponding to the map $F_2$ (see (\ref{affine})) are integers, and the locally defined map $\Pi_2\circ\theta^{-1}$ from $S^3\setminus I$ into $U_2$ extends to an elementary map ${\bf g}$ from $B^2\setminus I$ onto the Reinhardt domain $\Omega:=\Pi_2(\tilde T_1)$.

Let ${\cal L}_z':=\varphi^{-1}\left(B^2\cap{\cal L}_z\right)$, ${\cal L}_w':=\varphi^{-1}\left(B^2\cap{\cal L}_w\right)$, ${\cal L}_z'':=\varphi\left(B^2\cap{\cal L}_z\right)$, ${\cal L}_w'':=\varphi\left(B^2\cap{\cal L}_w\right)$. The map $\hat\varphi:={\bf g}\circ\varphi\circ\theta$ takes $\tilde T_1':=\tilde T_1\setminus\theta^{-1}\left(\Bigl({\cal L}_z'\cup{\cal L}_w'\Bigr)\setminus I\right)$ onto $\Omega\setminus{\bf g}\left(\Bigl({\cal L}_z''\cup{\cal L}_w''\Bigr)\setminus I\right)$. Recall that on an open subset of $\partial \tilde T_1'$ the map $\hat\varphi$ coincides with $f\circ\Pi_1$ and thus extends to all of $\tilde T_1$. Therefore, ${\bf g}$ extends to $B^2\cap I$, and $f\circ\Pi_1$ maps $\tilde T_1$ onto $\Bigl(\Omega\cup{\bf g}(B^2\cap I)\Bigr)\setminus{\bf g}\left({\cal L}_z''\cup{\cal L}_w''\right)$. Thus, $\hat D_2=\Biggl(\Bigl(\Omega\cup{\bf g}(B^2\cap I)\Bigr)\setminus{\bf g}\left({\cal L}_z''\cup{\cal L}_w''\right)\Biggr)\cup f(\hat D_1\cap I)$. Since $\hat D_2$ is a Reinhardt domain, it follows that ${\bf g}\Bigl(\left({\cal L}_z''\cup{\cal L}_w''\right)\setminus I\Bigr)\subset f(\hat D_1\cap I)$, and therefore $\hat D_2=\Omega\cup{\bf g}\left(B^2\cap I\right)$. In particular, $\hat D_1\cap I\ne\emptyset$.

Further, $\Omega$ is a bounded Reinhardt domain not intersecting $I$ such that there exists an elementary map from $B^2\setminus I$ onto $\Omega$. Therefore, $\Omega$ has the form (\ref{onemoreform}) with $\alpha_0=2\pi$, and ${\bf g}$ has the form (\ref{formmmmmmg}). 
If either $a_2>0$ and $c_2>0$, or $b_2>0$ and $d_2>0$, it follows from (\ref{formmmmmmg}) that ${\bf g}\left(B^2\cap I\right)=\{0\}$. However, $\Omega\cup\{0\}$ is not an open set in this case, and therefore either $a_2=0$, $d_2=0$, or $b_2=0$, $c_2=0$. Thus 
$$
\hat D_2=\left\{(z,w)\in\CC^2: C_2'|z|^{\frac{2}{a_2'}}+E_2' |w|^{\frac{2}{b_2'}}<1\right\}, 
$$
where $a_2',b_2'\in\NN$, $C_2'>0$, $E_2'>0$ (cf. (\ref{domm}) and (\ref{dommmmmm})), and ${\bf g}$ up to permutation of the variables has the form
$$
\begin{array}{lll}
z&\mapsto&\hbox{const}\, z^{a_2'},\\
w&\mapsto&\hbox{const}\, w^{b_2'}.
\end{array}
$$
This description shows that transformations in $G_2$ have the form (\ref{g2form}) with $\alpha_2'=2\pi/a_2'$, $\alpha_2=0$, $\beta_2'=0$, $\beta_2=2\pi/b_2'$.    

It is straightforward to observe that, since $\varphi\circ \hat G_1\circ \varphi^{-1}\subset \hat G_2$ and $\varphi(B^2\cap{\cal L}_z)\not\subset I$, $\varphi(B^2\cap{\cal L}_w)\not\subset I$, every transformation in $\hat G_1$ has the form
\begin{equation}
\begin{array}{lll}
z&\mapsto& e^{i\alpha}z,\\
w&\mapsto& e^{i\alpha}w,
\end{array}\label{alpha1eqalpha2}
\end{equation}
for $\alpha\in\RR$. Since $\hat D_1$ is a bounded Reinhardt domain intersecting $I$ with logarithmic diagram affinely equivalent to that of $B^2$, either up to permutation of the variables it has the form
\begin{equation}
\left\{(z,w)\in\CC^2: C_1|z|^{2a_1}|w|^{2c_1}+E_1 |w|^{2b_1}<1,\, w\ne 0\right\},\label{domainhatd1}
\end{equation}
or it has the form
\begin{equation}
\left\{(z,w)\in\CC^2: C_1|z|^{2a_1}+E_1 |w|^{2b_1}<1, \right\},\label{domainhatd12}
\end{equation} 
where $a_1,b_1, c_1\in\RR$, $a_1>0$, $b_1>0$, $c_1\le 0$ and $C_1>0$, $E_1>0$ (cf. (\ref{domainv})).

Assume first that $\hat D_1$ is a domain of the form (\ref{domainhatd1}). Then the group $G_1$ consists of transformations (\ref{groupg1forrrrm}) with $\alpha_1'=2\pi a_1$, $\alpha_1=0$, $\beta_1'=2\pi c_1$, $\beta_1=2\pi b_1$. Since all transformations in $\hat G_1$ are of the form (\ref{alpha1eqalpha2}), it follows that $a_1\in\NN$. Thus, the matrix $A$ corresponding to the map $F_1$ (see (\ref{affine})) up to permutation of the rows is
$$
\left(
\begin{array}{ll}
a_1 & c_1\\
0 & b_1
\end{array}
\right).
$$
For an appropriate choice of an element of $\theta^{-1}$, the locally defined map $\hat\varphi\circ F_1\circ\theta^{-1}$ from $\hat D_1\setminus I$ into $\hat D_2$ coincides with $f$, and hence extends to all of $\hat D_1$. It then follows from this representation of $f$ that either $f(\hat D_1\cap I)={\bf g}\left({\cal L}_z''\right)$ or $f(\hat D_1\cap I)={\bf g}\left({\cal L}_w''\right)$,  and therefore ${\bf g}\Bigl(\left({\cal L}_z''\cup{\cal L}_w''\right)\setminus I\Bigr) \not\subset f(\hat D_1\cap I)$, in contradiction to what we established above. This shows that in fact $\hat D_1$ has the form (\ref{domainhatd12}).

Since all transformations in $\hat G_1$ are of the form (\ref{alpha1eqalpha2}), it follows that $a_1,b_1\in\NN$. Hence the locally defined map $\theta\circ\Pi_1^{-1}$ from $\hat D_1\setminus I$ into $B^2\setminus I$ extends to an elementary map ${\bf h}$ from $\hat D_1$ onto $B^2$. Clearly, up to permutation of its components, the map ${\bf h}$ has the following form
$$
\begin{array}{lll}
z&\mapsto& \hbox{const}\, z^{a_1},\\
w&\mapsto& \hbox{const}\, w^{b_1},
\end{array}
$$
and we have $f={\bf g}\circ\varphi\circ{\bf h}$. It is straightforward to see that there exists no proper subdomain of $\hat D_1$ 
mapped properly by $f$ onto a bounded Reinhardt domain and whose envelope of holomorphy coincides with $\hat D_1$. Therefore, $D_1=\hat D_1$, and hence $D_2=\hat D_2$.

We thus have obtained (vi) of Theorem \ref{mainresult}. The proof of the theorem is now complete.\qed

{\obeylines
Department of Mathematics
The Australian National University
Canberra, ACT 0200
AUSTRALIA
E-mail address: Alexander.Isaev@maths.anu.edu.au
\hbox{ \ \ }
\hbox{ \ \ }
Department of Complex Analysis
Steklov Mathematical Institute
42 Vavilova St.
Moscow 117966
RUSSIA
E-mail address: kruzhil@mi.ras.ru
}

\end{document}